\documentclass[10pt]{art}

\usepackage[skip=0.4em]{parskip}
\usepackage[mathcal]{euscript}
\usepackage{pgfplots, caption, float, xargs, stmaryrd}

\newcommand{\Art}{\ensuremath{\mathsf{Art}}}
\newcommand{\Alg}{\ensuremath{\mathsf{Alg}}}
\newcommand{\B}{\ensuremath{\mathbf{B}}}
\newcommand{\A}{\ensuremath{\mathbf{A}}}
\newcommand{\K}{\ensuremath{\mathbb{K}}}

\newcommand{\M}{\ensuremath{\mathcal{M}}}

\newcommand{\PT}{\ensuremath{\mathsf{P}\mathbb{T}}}
\newcommand{\RT}{\ensuremath{\mathsf{R}\mathbb{T}}}
\newcommand{\Tbb}{\ensuremath{\mathbb{T}}}
\newcommand{\T}{\ensuremath{\mathcal{T}}}
\renewcommand{\vert}{\ensuremath{\textsf{vert}}}

\renewcommand{\P}{\ensuremath{\mathcal{P}}}
\renewcommand{\C}{\ensuremath{\mathcal{C}}}
\newcommand{\Pz}{\ensuremath{\mathcal{P}^{!}}}
\newcommand{\Pza}{\ensuremath{\mathcal{P}^\text{!`}}}

\newcommand{\Cza}{\ensuremath{\mathcal{C}^\text{!`}}}
\newcommand{\Pinf}{\ensuremath{\mathcal{P}_\infty}}
\renewcommand{\Bar}{\ensuremath{\mathbf{B}}}
\newcommand{\Cobar}{\ensuremath{\mathbf{\Omega}}}
\newcommand{\MC}{\ensuremath{\mathsf{MC}}}
\newcommand{\g}{\ensuremath{\mathfrak{g}}}
\newcommand{\h}{\ensuremath{\mathfrak{h}}}
\newcommand{\Sbb}{\ensuremath{\mathbb{S}}}
\newcommand{\Sh}{\ensuremath{\mathsf{Sh}}}
\renewcommand{\N}{\ensuremath{\mathbb{N}}}
\newcommand{\Nmod}{\ensuremath{\mathbb{N}\text{-mod}}}
\newcommand{\Smod}{\ensuremath{\mathbb{S}\text{-mod}}}
\newcommand{\Tw}{\ensuremath{\mathsf{Tw}}}
\newcommand{\Ass}{\ensuremath{\mathit{Ass}}}
\newcommand{\Com}{\ensuremath{\mathit{Com}}}
\newcommand{\Lie}{\ensuremath{\mathit{Lie}}}
\newcommand{\Linf}{\ensuremath{\mathsf{L}_\8}}
\newcommand{\Einf}{\ensuremath{\mathbb{E}_\8}}

\newcommand{\Cscr}{\ensuremath{\mathscr{C}}}

\renewcommand{\S}{\ensuremath{\mathcal{S}}}
\newcommand{\Sfin}{\ensuremath{\mathcal{S}_*^{\mathsf{fin}}}}
\newcommand{\Sp}{\ensuremath{\mathsf{Sp}}}
\newcommand{\Spf}{\ensuremath{\mathsf{Spf}}}

\newcommand{\Hoch}{\ensuremath{\mathsf{Hoch}}}
\newcommand{\Spec}{\ensuremath{\mathsf{Spec}}}

\newcommand{\FMP}{\ensuremath{\mathsf{FMP}}}
\newcommand{\FMPP}{\ensuremath{\FMP_\P}}
\newcommandx{\Tgt}[2][1]{\ensuremath{\mathbb{T}_{#1}(#2)}}
\newcommand{\Dfrak}{\ensuremath{\mathfrak{D}}}
\newcommand{\ArtP}{\ensuremath{\Art_\P}}

\title{Operads in Derived Deformation Theory}
\author{
  \textcolor{gray80}{Ramkumar Ramachandra} \\
  \itshape \textcolor{gray80}{Université Paris-Cité} \\
  \ttfamily \textcolor{gray80}{r@artagnon.com}
}
\date{}

\pgfplotsset{compat=1.17}
\DeclareCaptionFont{mdit}{\mdseries\itshape}
\color{gray65}

\begin{document}
\maketitle
\thispagestyle{empty}
\begin{abstract}
  A theorem by Pridham and Lurie provides an equivalence between formal moduli problems and Lie algebras in characteristic zero. In his work, Lurie has distilled the axioms that the algebras appearing in the formal moduli problem need to satisfy, and worked out the case of \Einf-algebras using an incarnation of the Koszul duality, in the setting of \8-operads. The more recent work of Calaque-Campos-Nuiten extends Lurie's work to obtain an equivalence between formal moduli problem parameterized by a colored operad, and algebras over its Koszul dual operad. This manuscript is both, a pedagogical exposition, and a questioning of their work, with modest, but original, supporting lemmas.
\end{abstract}

\tableofcontents
\newpage

\setcounter{section}{-1}
\renewcommand{\thesection}{\texorpdfstring{$\emptyset$}{0}}
\section{Introduction}
Derived deformation theory nests under a larger body of mathematics known as \emph{derived algebraic geometry}. In order to elaborate on the scope of our work, let us first provide a ``big picture'' view of derived geometry.

\subsection*{Derived geometry}
Derived geometry is a theory which generalizes the classical counterpart, while faithfully embedding the classical version at non-singular points in a space. We work over a \emph{moduli space}, which is, informally, a space that classifies something, say solutions to a set of equations, or birational classification of curves. Examples of moduli spaces include affine schemes, smooth manifolds, stacks, but perhaps the simplest example is $\mathbb{CP}^n$. The singularities in these spaces are of two different kinds.

\begin{enumerate}
  \item[(1)] \emph{Quotient singularities} appear in the quotient of a group action when a point has a stabilizer under the action.

    \begin{figure}[H]
      \centering
      \begin{tikzpicture}[scale=0.5]
        \begin{axis}[
            hide axis,
            ymin=-0.1,
            ymax=4.1
          ]
          \addplot[yaleblue, samples=200, domain=-2:2] {x^2};
          \addplot[yaleblue, samples=200, domain=-2:2] {x^2+1};
        \end{axis}
      \end{tikzpicture}
      \hskip 10pt
      \begin{tikzpicture}[scale=0.5]
        \begin{axis}[
            hide axis,
            ymin=-2,
            ymax=2
          ]
          \addplot[yaleblue, samples=200, domain=-2:2] {x^2};
          \addplot[yaleblue, samples=200, domain=-2:2] {-x^2};
          \addplot[raspeberry, samples=200, domain=-2:2] {0};
        \end{axis}
      \end{tikzpicture}
      \caption{The space on the left has zeroth order quotient singularity, while the space on the right has higher order quotient singularities.}
    \end{figure}
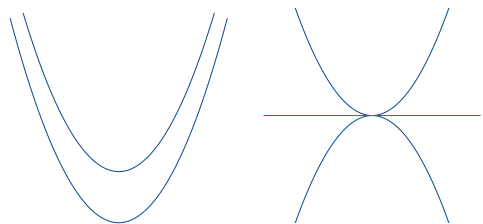

  \item[(2)] \emph{Intersection singularities} appear when intersections are not transverse.

    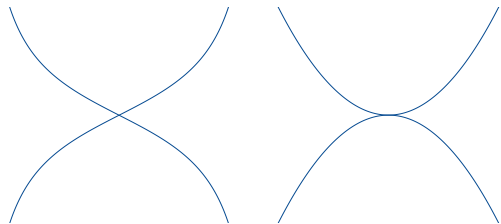
\begin{figure}[H]
      \centering
      \hskip 10pt
      \begin{tikzpicture}[scale=0.5]
        \begin{axis}[
            hide axis,
            ymin=-2.2,
            ymax=2.2
          ]
          \addplot[yaleblue, samples=200, domain=-1.5:1.5] {tan(deg(x))};
          \addplot[yaleblue, samples=200, domain=-1.5:1.5] {-tan(deg(x))};
        \end{axis}
      \end{tikzpicture}
      \begin{tikzpicture}[scale=0.5]
        \begin{axis}[
            hide axis,
            ymin=-2,
            ymax=2
          ]
          \addplot[yaleblue, samples=200, domain=-2:2] {x^2};
          \addplot[yaleblue, samples=200, domain=-2:2] {-x^2};
        \end{axis}
      \end{tikzpicture}
      \caption{The space on the left has zeroth order intersection singularity, which is to say that the intersection is transverse, while the space on the right has higher order intersection singularities.}
    \end{figure}
\end{enumerate}

A singularity can be characterized by saying that not all tangent vectors can be integrated into a curve. In the pursuit to study higher-order singularities, the tangent space was then enhanced to a \emph{tangent complex}, whose homology vanishes at non-singular points. At non-singular points, the structure of the tangent complex is trivial (Abelian), and at singular points, it encodes a lot of information about the singularity. It turns out that the tangent complex neatly breaks up into two parts: intersection singularities correspond to the positive part, and are studied in \emph{derived intersection theory} using \emph{derived schemes}, while the quotient singularities correspond to the negative part, and are studied in \emph{derived deformation theory} using \emph{derived stacks}.

The central role of \emph{Lie algebras} in this picture derives from the fact that for a tangent complex $\Tbb$, $\Tbb[-1]$ is always equipped with some kind of Lie algebra structure.

\subsection*{Subject of this manuscript}
The topic of this manuscript interests itself in the following central fact: the entire \emph{formal neighborhood} can be reconstructed from the tangent complex equipped with its Lie algebra structure. Research in the field today has led to a seminal result called the \emph{Lurie-Pridham equivalence}. Our work is based on a variation of this result by \emph{Calaque-Campos-Nuiten}~\cite{Calaque19} done in the setting of ordinary operads.

Derived deformation theory is a highly algebraic theory, and this has two consequences for us.

\begin{enumerate}
  \item[(a)] We will not concern ourselves with either derived schemes or stacks. Indeed, the only ``geometry'' that appears in our manuscript is the \8-category of spaces.
  \item[(b)] Much of the modern literature on the subject, and the work of \emph{Calaque-Campos-Nuiten} in particular, is written in the language of \emph{operads}. For this reason, we dedicate a significant portion of our manuscript exclusively to the study of operads.
\end{enumerate}

\subsection*{Motivation}
The objective of this manuscript is to understand the motivation behind the development of the machinery of \8-operads in the context of derived deformation theory.

\subsection*{Our contribution}
Our contribution is two-fold. For the more pedagogical contribution, we provide a firm grounding to understand \ref{sec:op-deformation}, by giving the reader a thorough understanding to the theory of operads in \ref{sec:operads} and a motivated storyline leading up to Lurie's results in \ref{sec:derived-deformation}.

We present the less pedagogical and more original contribution in \ref{sec:op-deformation}, where we fill the reader in on several details, and obtain modest supporting results.

\subsection*{Prerequisites}
The audience for this graduate-level manuscript is expected to have a fair bit of mathematical background, up to an understanding of \8-categories~\cite{Lurie09}. No exposure to the theory of operads is assumed.

\subsection*{Acknowledgements}
The author would like to thank Ricardo Campos~\footnote{CNRS Researcher, Université Paul Sabatier, Toulouse} for his continued guidance and support.

\newpage
\renewcommand{\thesection}{\Roman{section}}
\section{Operads\label{sec:operads}}
An algebra can be characterized by an abstract multiplication operator $\mu$, that satisfies certain properties. For example, tensor algebra is characterized by the properties of the tensor product, of either vector spaces or modules. We generalize this statement by saying that vector spaces or modules form a symmetric monoidal category with respect to the tensor product, $\otimes$. The theory of \emph{operads} is a further generalization, whence once considers multiplication operators with multiple inputs and one output. Once sufficient machinery has been built, we can discuss algebras that are ``homotopy coherent'', or \Pinf-algebras.

\begin{caution}[$1$-operadic theory versus \8-operadic theory]
  \emph{Colored operads} consider operads with multiple inputs and multiple outputs: the ordinary theory of operads nests within this theory, as $1$-colored operads. Finally, \8-operads generalize colored operads, to yield a theory in which we can define an \Einf-algebra, and do non-commutative derived deformation theory.

  The theory of \8-operads is developed in detail in \cite{Lurie12}, and involves a fair bit of machinery. Fortunately, the theory of $1$-operads is sufficient for our purposes, and we provide an exposition to this theory, concluding with \Pinf-algebras. Even though the notation might suggest otherwise, \Pinf-algebras and \8-quasi-isomorphisms are notions of $1$-operadic theory.
\end{caution}

\begin{convention}\label{con:char0k}
  We will always assume the ground field, $\K$, to be of characteristic zero. Unless indicated otherwise, tensor products will be over $\K$.
\end{convention}

\begin{convention}\label{con:01triv}
  All operads $\P$, and cooperads $\C$ will be assumed to be trivial in degrees $0$ and $1$: $\P_0 = 0$, $\P_1 = \I$, $\C_0 = 0$, and $\C_1 = \I$, where $\I$ denotes the unit of the underlying monoidal category.
\end{convention}

\begin{intuition}[A high-level summary]
  Given operad $\P$ and its \emph{Koszul dual} \Pz, the following adjunction

  \begin{equation*}
    \begin{tikzcd}
      \Bar : \P\textsf{-algebras} \arrow[r, bend left] \arrow[r, phantom, "\bot" description] & \textsf{conilpotent }\Pz\textsf{-coalgebras} \arrow[l, bend left] : \Cobar
    \end{tikzcd}
  \end{equation*}

  is a Quillen equivalence. Let $A$ be a \P-algebra. Then,

  \begin{align*}
    \Bar A                & = \Pza \circ A                 \\
    \textsf{where }\Pz(n) & := \Hom_\K(\Pza(n), \K)[1 - n]
  \end{align*}

  Now,

  \begin{equation*}
    \begin{matrix}
       & \Map(\Cobar\Pza, \End{A})                   & = & \MC(\g)                &       \\
       & \cap                                        &   & \cap                   &       \\
       & \Map_{\mathsf{non-dg}}(\Cobar\Pza, \End{A}) & = & \Hom_\K(\Pza, \End{A}) & := \g
    \end{matrix}
  \end{equation*}

  where $\MC(\g)$ are the \emph{Maurer-Cartan elements} of \g. Moreover, there exists a operad morphism, which induces "restriction of scalars" functor at the level of algebras:

  \begin{equation*}
    \Linf \rightarrow \P \otimes \Pz
  \end{equation*}

  where $\Linf$ is the cofibrant resolution of the Lie operad.
\end{intuition}

\subsection{Non-symmetric operads}
An algebraic operad $\P$ is a functor from $\Vect$ to \Vect, that sends algebras over the source vector space to \emph{free algebras} over the destination vector space. $\P$ consists of the data of the triple $(\P, \gamma, \eta)$, where $\gamma: \P \circ \P \rightarrow \P$ and $\eta: \I \rightarrow \P$, forming a \emph{monoid} [\ref{def:monoid}] in a strict monoidal category. In particular, a monoid in a strict monoidal category of endofunctors is called a \emph{monad} [\ref{def:monad}].

Non-symmetric operads, or \emph{ns-operads}, are a simple setting we can use to initiate the study of operads. This setting will then allow us to study the theory of symmetric operads, which we otherwise just call operads.

\begin{definition}[Schur functor]
  Let $M_\bullet = \{M_n\}_{n \geq 0}$ be an arity-graded vector space. The \emph{Schur functor} $M: \Vect \rightarrow \Vect$ associated to $M_\bullet$, is defined as:

  \begin{equation*}
    \M(V) := \bigoplus_n \M_n \otimes V^{\otimes n}
  \end{equation*}
\end{definition}

\begin{definition}[Non-symmetric operad\label{def:nsop}]
  An ns-operad $(\P, \gamma, \alpha)$ is defined as a Schur functor that is compatible with the monoidal structure

  \begin{equation*}
    \P(n) := \bigoplus_n \P_n \otimes V^{\otimes n}
  \end{equation*}

  In other words, an ns-operad, $(P, \gamma, \alpha)$, is an arity-graded vector space $\P = \{\P_n\}_{n \geq 0}$ equipped with the composition maps

  \begin{equation*}
    \gamma_{i_1, \ldots, i_k} : \P_k \otimes \P_{i_1} \otimes \ldots \otimes \P_{i_k} \rightarrow \P_{i_1 + \ldots + i_k}
  \end{equation*}

  and an element $\id \in \P_1$, such that the transformation functors $\gamma: \P \circ \P \rightarrow \P$ and $\gamma: \I \rightarrow \P$ make $(\P, \gamma, \alpha)$ into a monoid.
\end{definition}

We use the formalism of planar rooted trees to give an alternative definition of an ns-operad.

\begin{definition}[Combinatorial definition of ns-operad]
  The category of graded vector spaces is denoted by \Nmod. Let $M$ be an arity-graded vector space with $M_0 = 0$. For any \emph{planar rooted tree} [\ref{sec:trees}] $t$,

  \begin{equation*}
    M_t := \bigotimes_{v \in \vert(t)} M_{|v|}
  \end{equation*}

  We then get a functor

  \begin{align*}
    \PT      & : \Nmod \rightarrow \Nmod      \\
    \PT(M)_n & := \bigoplus_{t \in \PT_n} M_t \\
  \end{align*}

  The corolla enables us to define $\eta: \I_{\Nmod} \rightarrow \PT$, and the \emph{substitution of trees}, given by grafting trees with matching inputs, enables us to define $\alpha: \PT \circ \PT \rightarrow \PT$.

  Then, the definition of an ns-operad is an arity-graded module \P, together with a map $\PT(\P) \rightarrow \P$ compatible with $\alpha$ and $\eta$ in the usual sense.
\end{definition}

\subsection{Symmetric operads}
\begin{definition}[Non-dg \Sbb-module\label{def:nondg-smod}]
  A non-dg \Sbb-module is defined as the family

  \begin{equation*}
    M := \{M(0), M(1), \ldots, M(n), \ldots\}
  \end{equation*}

  where each $M(n)$ is a right \K[$\Sbb_n$] module. Since, $\K$ is a field [\ref{con:char0k}], an \Sbb-module is an arity-graded vector space.
\end{definition}

\begin{definition}[Schur functor]
  The Schur functor is an endofunctor over \Vect

  \begin{align*}
    \widetilde{M}    & : \Vect \rightarrow \Vect                          \\
    \widetilde{M(V)} & := \bigoplus_n M(n) \otimes_{\Sbb_n} V^{\otimes n}
  \end{align*}
\end{definition}

\begin{lemma}[Constructions on the Schur functor]
  The direct sum, tensor product, and composition of two Schur functors $F, G: \Vect \rightarrow \Vect$ are given by

  \begin{align*}
    (F \oplus G)(V)  & := F(V) \oplus G(V)  \\
    (F \otimes G)(V) & := F(V) \otimes G(V) \\
    (F \circ G)(V)   & := F(V) \circ G(V)
  \end{align*}
\end{lemma}

\begin{lemma}[Direct sum of non-dg \Sbb-modules]
  \begin{equation*}
    (M \oplus N)(n) := M(n) \oplus N(n)
  \end{equation*}
\end{lemma}

\begin{lemma}[Tensor product of non-dg \Sbb-modules]
  \begin{equation*}
    (M \otimes N)(n) := \bigoplus_{i + j = n} M_i \otimes N_j \otimes \K[\Sh(i, j)]
  \end{equation*}

  where $\Sh(i, j)$ is the $ij$-\emph{shuffle} [\ref{def:shuffle}].
\end{lemma}

\begin{lemma}[Composition of non-dg \Sbb-modules]
  \begin{equation*}
    (M \circ N)(n) := \bigoplus_{k \geq 0} M_k \otimes_{\Sbb_k} N^{\otimes k}
  \end{equation*}
\end{lemma}

Compare the above definition with \ref{def:nsop}. Replacing $\P_n$ with an \emph{\Sbb-module}, we get the general definition of an operad.

\begin{definition}[Algebraic operad, re-stated]
  An algebraic operad, \P, consisting of the data of the triple $(\P, \gamma, \eta)$, is a Schur functor that is compatible with monad laws.

  \begin{align*}
    \P(V) := & \bigoplus_n \P_n \otimes V^{\otimes n}           \\
    =        & \bigoplus_n \P(n) \otimes_{\Sbb_n} V^{\otimes n}
  \end{align*}
\end{definition}

\begin{notation}[$\P_n$]
  We use two different notations to clarify the framework we're working in. In the case of symmetric operads,

  \begin{equation*}
    \P(n) := \P_n \otimes \K[\Sbb_n]
  \end{equation*}
\end{notation}

\begin{definition}[The tree monad]
  Let $M$ be an \Sbb-module. We define $\T: \Vect \rightarrow \Vect$ as

  \begin{align*}
    \T_0 M & := \I                                \\
    \T_1 M & := \I \oplus M                       \\
    \T_2 M & := \I \oplus (M \circ (\I \oplus M)) \\
    \ldots                                        \\
    \T_n M & := \I \oplus (M \circ \T_{n - 1} M)
  \end{align*}

  By definition, the tree module $\T M$ over an \Sbb-module is

  \begin{equation*}
    \T M := \colim_n \T_n M
  \end{equation*}

  There is an operad structure $\gamma$ on $\T M$ such that $\T(M) := (\T M, \gamma, j)$, where $j: M \rightarrow \T M$ sends $M \circ \T_{n - 1} M \mapsto \T_n M$, and $\gamma: \T M \rightarrow \T M$. $\T$ is hence a \emph{free operad} on $M$.
\end{definition}

\subsection{Modern definition of an operad}
To summarize everything we've covered so far, let us now present the abstract definition of an operad, from two different equivalent viewpoints. In particular, we introduce grafting of trees, which will be used to generalize operads to the differential graded setting.

\begin{definition}[Combinatorial definition of an operad]
  Let $M$ be an arity-graded vector space with $M_0 = 0$. For any \emph{rooted tree} [\ref{sec:trees}] $t$, we define the \emph{treewise tensor product} $M_t$ as

  \begin{equation*}
    M_t := \bigotimes_{v \in \vert(t)} M_{|v|}
  \end{equation*}

  and a functor

  \begin{align*}
    \Tbb       & : \Smod \rightarrow \Smod        \\
    \Tbb(M)(X) & := \bigoplus_{t \in \RT(X)} M(t)
  \end{align*}

  We then construct a monad structure on $\Tbb$ by \emph{substitution of trees}. In order to do this, we need to define $\iota: \id_{\Smod} \rightarrow \Tbb$ and $\alpha: \Tbb \circ \Tbb \rightarrow \Tbb$ to obtain $(\Tbb, \iota, \alpha)$.

  Given \Sbb-module $M$, $\iota$ sends $M(X) \mapsto \Tbb(M)(X)$. In $\RT(X)$, we have one particular tree corresponding to the corolla, so that $M(\mathsf{corolla}) = M(X) = \bigoplus \Tbb(M)(X)$.

  Next, let us illustrate the substitution of trees that defines $\alpha$: we graft matching inputs onto the leaves.

  
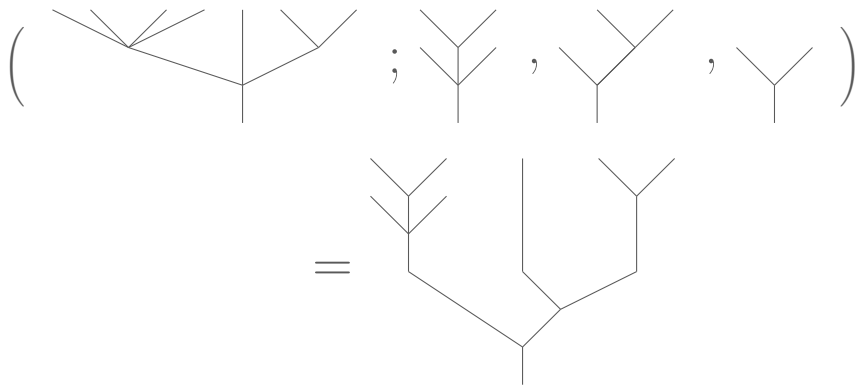
\begin{figure}[H]
    \centering
    \begin{tikzpicture}[scale=0.5]
      \node at (-1, 1.5) {$\Biggl($};
      \draw (5, 0) to (5, 1);
      \draw (5, 1) to (2, 2);
      \draw (5, 1) to (5, 3);
      \draw (5, 1) to (7, 2);
      \draw (2, 2) to (0, 3);
      \draw (2, 2) to (1, 3);
      \draw (2, 2) to (3, 3);
      \draw (2, 2) to (4, 3);
      \draw (7, 2) to (6, 3);
      \draw (7, 2) to (8, 3);
      \node[font=\huge] at (9, 1.5) {;};
    \end{tikzpicture}
    \begin{tikzpicture}[scale=0.5]
      \draw (1, 0) to (1, 1);
      \draw (1, 1) to (1, 2);
      \draw (1, 1) to (0, 2);
      \draw (1, 1) to (2, 2);
      \draw (1, 2) to (0, 3);
      \draw (1, 2) to (2, 3);
      \node[font=\huge] at (3, 1.5) {,};
    \end{tikzpicture}
    \begin{tikzpicture}[scale=0.5]
      \draw (1, 0) to (1, 1);
      \draw (1, 1) to (0, 2);
      \draw (1, 1) to (2, 2);
      \draw (1, 1) to (2, 2);
      \draw (2, 2) to (1, 3);
      \draw (2, 2) to (3, 3);
      \node[font=\huge] at (4, 1.5) {,};
    \end{tikzpicture}
    \begin{tikzpicture}[scale=0.5]
      \draw (1, 0) to (1, 1);
      \draw (1, 1) to (0, 2);
      \draw (1, 1) to (2, 2);
      \node at (3, 1.5) {$\Biggr)$};
    \end{tikzpicture}
    \newline\newline
    \begin{tikzpicture}[scale=0.5]
      \node[font=\huge] at (-1, 3) {=};
      \draw (4, 0) to (4, 1);
      \draw (4, 1) to (1, 3);
      \draw (4, 3) to (5, 2);
      \draw (4, 1) to (5, 2);
      \draw (5, 2) to (7, 3);
      \draw (1, 3) to (1, 4);
      \draw (4, 3) to (4, 6);
      \draw (1, 4) to (1, 5);
      \draw (1, 4) to (0, 5);
      \draw (1, 4) to (2, 5);
      \draw (1, 5) to (0, 6);
      \draw (1, 5) to (2, 6);
      \draw (7, 3) to (7, 5);
      \draw (7, 5) to (6, 6);
      \draw (7, 5) to (8, 6);
    \end{tikzpicture}
  \end{figure}
\end{definition}

\begin{definition}[Partial definition of an operad]
  Let $\P$ be an operad and $\mu \in \P(n)$, $\nu \in \P(m)$ be two operations. We then define the partial composite $\circ_i$ of $\mu$ and $\nu$ as

  \begin{align*}
    - \circ_i -     & : \P(m) \otimes \P(n) \rightarrow \P(m + n - 1)         \\
    \mu \circ_i \nu & := \gamma(\mu; \id, \ldots, \id, \nu, \id, \ldots, \id)
  \end{align*}

  This can be represented as a grafting of trees, where the root of $\nu$ is grafted on to the $i$th leaf of $\mu$.

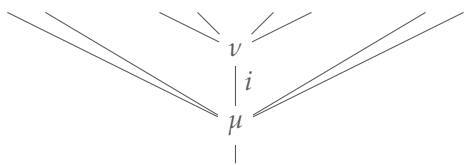
  \begin{figure}[H]
    \centering
    \begin{tikzpicture}[scale=0.5]
      \node (mu) at (6, 1) {$\mu$};
      \node (nu) at (6, 3) {$\nu$};
      \node[right] at (6, 2.2) {$i$};
      \draw (mu) to (nu);
      \draw (6, 0) to (mu);
      \draw (0, 4) to (mu);
      \draw (1, 4) to (mu);
      \draw (4, 4) to (nu);
      \draw (5, 4) to (nu);
      \draw (7, 4) to (nu);
      \draw (8, 4) to (nu);
      \draw (11, 4) to (mu);
      \draw (12, 4) to (mu);
    \end{tikzpicture}
  \end{figure}

  In order to describe the action of the symmetric groups and complete the definition of $\circ_i$, we axiomatize the two different cases of two-stage compositions.

  \begin{equation*}
    \begin{matrix}
      \textbf{(I)}  & (\lambda \circ_i \mu) \circ_{i + j - 1} \nu & = \lambda \circ_i (\mu \circ_j \nu) & 1 \leq i \leq l, 1 \leq j \leq m \\
      \textbf{(II)} & (\lambda \circ_i \mu) \circ_{k + m - 1} \nu & = (\lambda \circ_k \nu) \circ_i \mu & 1 \leq i \leq k \leq l           \\
    \end{matrix}
  \end{equation*}

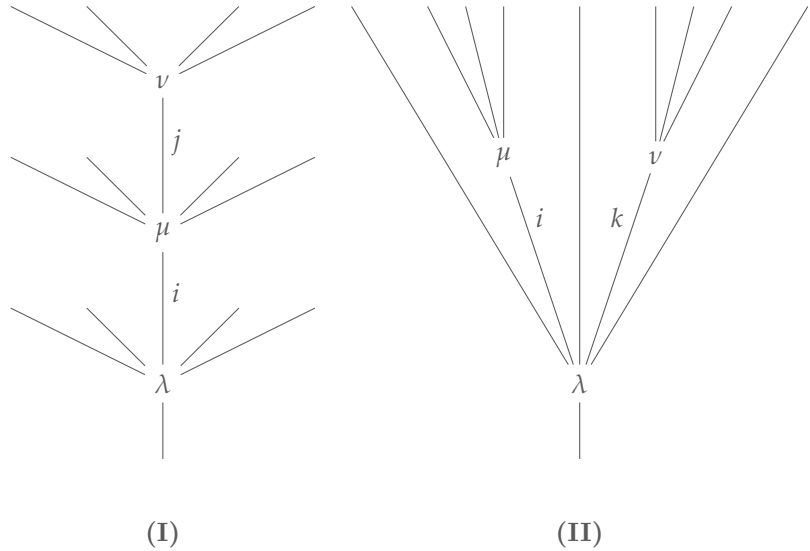
  \begin{figure}[H]
    \centering
    \begin{tikzpicture}
      \node (mu) at (2, 1) {$\lambda$};
      \node (nu) at (2, 3) {$\mu$};
      \node (lambda) at (2, 5) {$\nu$};
      \node[right] at (2, 2.2) {$i$};
      \node[right] at (2, 4.2) {$j$};
      \draw (mu) to (nu);
      \draw (nu) to (lambda);
      \draw (2, 0) to (mu);
      \draw (0, 2) to (mu);
      \draw (1, 2) to (mu);
      \draw (3, 2) to (mu);
      \draw (4, 2) to (mu);
      \draw (0, 4) to (nu);
      \draw (1, 4) to (nu);
      \draw (3, 4) to (nu);
      \draw (4, 4) to (nu);
      \draw (0, 6) to (lambda);
      \draw (1, 6) to (lambda);
      \draw (3, 6) to (lambda);
      \draw (4, 6) to (lambda);
      \node at (2, -1) {\textbf{(I)}};
    \end{tikzpicture}
    \hskip 10pt
    \begin{tikzpicture}
      \node (lambda) at (2, 1) {$\lambda$};
      \node (mu) at (1, 4) {$\mu$};
      \node (nu) at (3, 4) {$\nu$};
      \node[right] at (1.3, 3.2){$i$};
      \node[left] at (2.7, 3.2) {$k$};
      \draw (lambda) to (mu);
      \draw (lambda) to (nu);
      \draw (2, 0) to (lambda) to (2, 6);
      \draw (-1, 6) to (lambda);
      \draw (5, 6) to (lambda);
      \draw (0, 6) to (mu);
      \draw (0.5, 6) to (mu);
      \draw (1, 6) to (mu);
      \draw (4, 6) to (nu);
      \draw (3.5, 6) to (nu);
      \draw (3, 6) to (nu);
      \node at (2, -1) {\textbf{(II)}};
    \end{tikzpicture}
  \end{figure}

  The partial definition of an operad is then an \Sbb-module $\P$ equipped with partial compositions $\circ_i$ satisfying the compatibility relations (I) and (II) described above. It is also assumed that there is an element $\id \in \P(1)$ satisfying $\id \circ_i \nu = \nu$ and $\mu \circ_i \id = \mu$.
\end{definition}

\subsection{Algebra over operad}
\begin{definition}[\P-algebra]
  An algebra over operad \P, or a \P-algebra, is defined as a vector space $A$ equipped with map $\gamma_A: \P(A) \rightarrow A$ such that the following diagrams commute:

  
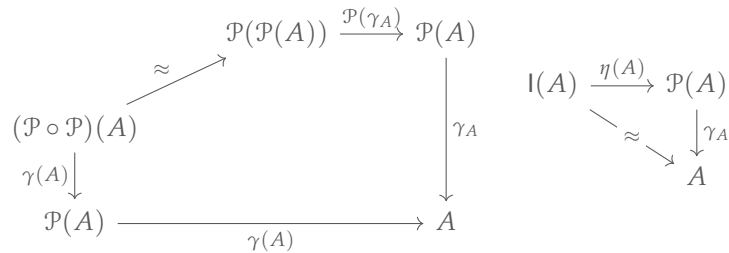
\begin{equation*}
    \begin{tikzcd}
      & \P(\P(A)) \arrow[r, "{\P(\gamma_A)}"] & \P(A) \arrow[dd, "{\gamma_A}"] \\
      (\P \circ \P)(A) \arrow[ur, "{\approx}"] \arrow[d, "{\gamma(A)}"{left}] & & \\
      \P(A) \arrow[rr, "{\gamma(A)}"{below}] & & A \\
    \end{tikzcd}
    \hskip 10pt
    \begin{tikzcd}
      \I(A) \arrow[r, "{\eta(A)}"] \arrow[dr, "{\approx}"{description}] & \P(A) \arrow[d, "{\gamma_A}"] \\
      & A \\
    \end{tikzcd}
  \end{equation*}

  See also \ref{def:monad-algebra} for the general definition of a monad algebra.
\end{definition}

\begin{definition}[Endomorphism operad]
  For any vector space $V$, the endomorphism operad is given by:

  \begin{equation*}
    \End{V}(n) := \Hom(V^{\otimes n}, V)
  \end{equation*}

  By convention, $V^{\otimes 0} = \K$. The right action of $\Sbb_n$ on $\End{V}(n)$ is induced by left action on $V^{\otimes n}$. The composition $\gamma$ is given by composition of endomorphisms:

  \begin{equation*}
    \begin{tikzcd}
      V^{\otimes i_1} \otimes \arrow[d, "{f_1}"] & \ldots & \otimes V^{\otimes i_k} \arrow[d, "{f_k}"] & = & V^{\otimes n} \arrow[d, "{f_1 \otimes \ldots \otimes f_k}"]\\
      V \otimes & \ldots \arrow[d, "{f}"] & \otimes V & = & V^{\otimes k} \arrow[d, "{f}"]\\
      & V & & = & V \\
    \end{tikzcd}
  \end{equation*}

  Evidently $\End{V}$ is an operad.
\end{definition}

\begin{lemma}[Alternative definition of \P-algebra]
  The \P-algebra structure on vector space $V$ is equivalent to the morphism of operads $\P \rightarrow \End{V}$.
\end{lemma}

We now supply three examples of \P-algebras, namely \Ass, \Com, and \Lie.

\begin{example}[\Ass]
  Let $\Ass: \Vect \rightarrow \Vect$ given by $\Ass(V) := \overline{T}(V)= \oplus_{n \geq 1} V^{\otimes n}$ (reduced tensor algebra). As an \Sbb-module, we have $\Ass(n) = \K[\Sbb_n]$, since $\K[\Sbb_n] \otimes_{\Sbb_n} V^{\otimes n} = V^{\otimes n}$ for $n \geq 1$, and $\Ass(0) = 0$. The map $\gamma(V): \Ass(\Ass(V)) \rightarrow \Ass(V)$ is given by ``composition of non-commutative polynomials''. This is a symmetric operad encoding associative algebras, since an algebra over $\Ass$ is the mapping $\overline{T}(A) \rightarrow A$.
\end{example}

\begin{example}[\Com]
  Let $\Com: \Vect \rightarrow \Vect$ be a Schur functor given by

  \begin{equation*}
    \Com(V) := \overline{S}(V) = \bigoplus_{n \geq 1} S^n V = \bigoplus_{n \geq 1} (V^{\otimes n})_{\Sbb_n}
  \end{equation*}

  As an \Sbb-module, we have $\Com(n) = \K$ with trivial action since

  \begin{equation*}
    \K \otimes_{\Sbb_n} V^{\otimes n} = S^n V \text{, for $n \geq 1$}
  \end{equation*}

  and $\Com(0) = 0$. The map $\gamma(V): \Com(\Com(V)) \rightarrow \Com(V)$ is given by ``composition of polynomials''. This is a symmetric operad encoding commutative algebras, since an algebra over $\Com$ is nonunital commutative.
\end{example}

\begin{example}[\Lie]
  Let $\Lie: \Vect \rightarrow \Vect$ be a Schur functor such that the space $\Lie(V) \subset \overline{T}(V)$ is generated by $V$ under the bracketing operation $[x, y] := xy - yx$. This is the free algebra over $V$. Let $\Lie(n)$ be the multilinear part of degree $n$ in free Lie algebra $\Lie(\K x_1 \oplus \ldots \oplus \K x_n)$. It can be shown that there is an operad structure on Schur functor $\Lie$ induced by the operad structure on $\Ass$. An algebra over $\Lie$ is a Lie algebra.
\end{example}

\subsection{The differential graded setting}
In this section, we generalize the notions of non-dg \Sbb-modules and operads to include a non-trivial differential. This allows us to continue discussing operads in full generality, in the following sections.

\begin{notation}[Infinitesimal composite of \Sbb-modules]
  Drawing inspiration from the partial definition of operads, in which $\circ_i$ is defined, we generalize this notation to define $\circ_{(1)}$, denoting the grafting of trees at one unspecified leaf. More precisely, given \Sbb-modules $M$ and $N$, we define their infinitesimal composite, $M \circ_{(1)} N$ as

  \begin{equation*}
    M \circ_{(1)} N := M \circ (\I; N)
  \end{equation*}

  This can be shown diagrammatically as

  
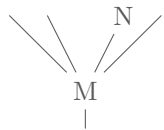
\begin{figure}[H]
    \centering
    \begin{tikzpicture}[scale=0.5]
      \node (M) at (2, 1) {M};
      \node (N) at (3, 3) {N};
      \draw (2, 0) to (M);
      \draw (M) to (0, 3);
      \draw (M) to (1, 3);
      \draw (M) to (N);
      \draw (M) to (4, 3);
    \end{tikzpicture}
  \end{figure}

  Elements of $M \circ_{(1)} N$ are of the form $(\mu; \id, \ldots, \id, \nu, \id, \ldots, \id)$
\end{notation}

\begin{notation}[Infinitesimal composite of morphisms]
  Given two \Sbb-module morphisms $f : M_1 \rightarrow M_2$ and $g : N_1 \rightarrow N_2$, we define $f \circ' g$ to be

  \begin{align*}
    f \circ' g & : M_1 \circ N_1 \rightarrow M_2 \circ (N_1; N_2)                                                                      \\
    f \circ' g & := \sum_i f \otimes (\id_{N_1} \otimes \ldots \otimes \id_{N_1} \otimes g \otimes \id_{N_1} \otimes \ldots \id_{N_1})
  \end{align*}

  This is unambiguously defined for a given $f$ and $g$: $g$ must occupy the $i$th position in the tensor product.
\end{notation}

\begin{definition}[Differential graded \Sbb-module]
  A \emph{graded \Sbb-module} $M$ is an \Sbb-module in the category of graded vector spaces. The arity $n$ is a graded $\Sbb_n$-module $\{M_p(n)\}_{p \in \Z}$. Equivalently, $M$ can be considered a family of \Sbb-modules. By abuse of notation,

  \begin{equation*}
    M := M_\bullet = \ldots \oplus M_0 \oplus \ldots \oplus M_p \oplus \ldots
  \end{equation*}

  A \emph{differential graded \Sbb-module} $(M, d)$ is a graded \Sbb-module $M$ equipped with differential $d$ of $\Sbb_n$-modules, such that $d^2 = 0$. Note that $H_\bullet(M)$ of a differential graded \Sbb-module form a graded \Sbb-module.

  The composite of two differential graded \Sbb-modules $(M, d_M)$ and $(N, d_N)$ comes with a differential $d_{M \circ N}$, which can be written using the notations we have introduced above.

  \begin{equation*}
    d_{M \circ N} := d_M \circ \id_N + \id_M \circ' d_N
  \end{equation*}
\end{definition}

It is easily checked that the category $(\textsf{dg } \Smod, \circ, \I)$ forms a monoidal category. We will now proceed to define a dg algebraic operad \P, and a dg \P-algebra.

\begin{definition}[Differential graded operad]
  A \emph{differential graded operad}, denoted $(\P, \gamma, \eta)$, is a monoid in the monoidal category $(\textsf{dg } \Smod, \circ, \I)$ where $\gamma : \P \circ \P \rightarrow \P$ is the composite map and $\eta : \I \rightarrow \P$ is the unit map, both of degree $0$, such that the following diagram commutes.

  \begin{equation*}
    \begin{tikzcd}
      \P \circ \P \arrow[r, "{\gamma}"] \arrow[d, "{d_{\P \circ \P}}"{left}] & \P \arrow[d, "{d_\P}"] \\
      \P \circ \P \arrow[r, "{\gamma}"] & \P \\
    \end{tikzcd}
  \end{equation*}

  where the differential $d_\P$ is given by

  \begin{equation*}
    \begin{matrix}
      d_\P(\gamma(\mu; \mu_1, \ldots, \mu_k)) := & \gamma(d_\P(\mu); \mu_1 \ldots, \mu_k) +                                                  \\
                                                 & \sum_{i = 1}^{k} (-1)^{\epsilon_i} \gamma(\mu; \mu_1, \ldots, d_\P(\mu_i), \ldots, \mu_k) \\
    \end{matrix}
  \end{equation*}

  and $\epsilon_i$ is defined as

  \begin{equation*}
    \epsilon_i := |\mu| + |\mu_1| + \ldots + |\mu_{i - 1}|
  \end{equation*}

  Since this is a chain complex, ${d_\P}^2 = 0$. Further, it can be verified that $H_\bullet(\P)$ carries a natural operadic structure.
\end{definition}

\begin{remark}
  The category $(\P, \gamma, \eta)$ forms a monad, just like the non-dg counterpart.
\end{remark}

\begin{definition}[Differential graded cooperad]
  A \emph{differential graded cooperad}, denoted $(\C, \Delta, \epsilon)$, is a comonoid in the monoidal category $(\textsf{dg } \Smod, \circ, \I)$ where $\Delta : \C \rightarrow \C \circ \C$ is the \emph{decomposition map} and $\epsilon : \C \rightarrow \I$ is the \emph{counit map}, both of degree $0$, such that the following diagram commutes.

  \begin{equation*}
    \begin{tikzcd}
      \C \arrow[r, "{\Delta}"]\arrow[d, "{d_\C}"{left}] & \C \circ \C \arrow[d, "{d_{\C \circ \C}}"] \\
      \C \arrow[r, "{\Delta}"] & \C \circ \C \\
    \end{tikzcd}
  \end{equation*}

  The formula for $\Delta(d_\C)$ is straightforward to work out.
\end{definition}

\begin{terminology}[Coaugmented cooperad]
  A cooperad is termed \emph{coaugmented} when it is endowed with a map $\eta: \I \rightarrow \C$ of degree $0$.
\end{terminology}

\begin{intuition}[Conilpotent cooperad]
  Conilpotency is a general condition on a cooperad that will be necessary to state certain facts in the following sections. However, it is strictly weaker than the convention we've adopted [\ref{con:01triv}], whereby $\C_0$ and $\C_1$ are trivial, and we ignore the qualifier altogether.
\end{intuition}

\begin{summary}
  This concludes our section detailing the machinery required to define the dg counterpart to operads. Since we have shown that the differential-graded version is an strict generalization of the non-dg version, and will proceed to drop the ``dg'' qualifiers everywhere.
\end{summary}

\subsection{Bar and cobar construction}
In classical homological algebra, there are several ways to obtain cofibrant and fibrant resolutions of modules: projective, injective and free resolutions make an early appearance. In operadic homological algebra, we use one canonical resolution called the \emph{Bar-cobar resolution}, which serves as the cofibrant resolution of the operadic algebra. We will not concern ourselves with fibrant resolutions or resolutions of operadic coalgebras, for reasons that will eventually become evident.

\begin{definition}[Convolution operad]
  Given operad $\P$ and cooperad $\C$, we make the following observations.

  \begin{align*}
    (\Hom_\Sbb(\C, \P), \star, \partial) & \;\;\text{forms a pre-Lie algebra} \\
    (\Hom_\Sbb(\C, \P), [,], \partial)   & \;\;\text{forms a Lie algebra}
  \end{align*}

  Here, $\Hom_\Sbb(\C, \P)$ is called a \emph{convolution algebra}, $\star$ is the pre-Lie multiplication operation, and $[,]$ is the Lie bracket. See \cite{Loday12} Section 6.4 for more.
\end{definition}

\begin{terminology}[Twisting morphism]
  The Maurer-Cartan equation in $\Hom_\Sbb(\C, \P)$ reads as either

  \begin{align*}
    \partial(\alpha) + \alpha \star \alpha = 0          \\
    \partial(\alpha) + \frac{1}{2} [\alpha, \alpha] = 0 \\
  \end{align*}

  depending on whether $\C$ and $\P$ are pre-Lie algebras or Lie algebras.

  A solution $\alpha: \C \rightarrow \P$ of degree $+1$ (in cohomological convention) is termed an \emph{operadic twisting morphism}, and is written $\Tw(\C, \P)$. See also: \ref{def:mc-elem}.
\end{terminology}

\begin{discussion}
  If $\alpha$ is of degree $-a$, consider the degree of the various terms in the Maurer-Cartan equation: $\partial(\alpha)$ is of degree $a + 1$, and $\frac{1}{2} [\alpha, \alpha]$ is of degree $2a$. We hence see that it is natural to define $\Tw$ to be $(+1)$-degree solutions of the Maurer-Cartan equation.
\end{discussion}

\begin{definition}[Minimal model of operad]
  Given operad $\P$, an operad $\M$ that is quasi-isomorphic to $\P$ is termed a  \emph{model} of operad $\P$. A given $\M = \T(E)$ is termed \emph{minimal} when $d(E) \subset \T(E)^{\geq 2}$.
\end{definition}

\begin{definition}[Twisted convolution operad\label{def:minmodel-op}]
  In any dg-Lie algebra, a $(+1)$-degree solution of the Maurer-Cartan equation gives rise to a twisted differential $\partial_\alpha$ on $\Hom_\Sbb(\C, \P)$. We can then define

  \begin{align*}
    \partial_\alpha(f)       & := \partial(f) + [\alpha, f]                   \\
    \Hom_\Sbb^\alpha(\C, \P) & := (\Hom_\Sbb^\alpha(\C, \P), \partial_\alpha)
  \end{align*}
\end{definition}

\begin{definition}[$\Bar$ and \Cobar]
  The functors $\Bar$ and $\Cobar$ operate as follows.

  \begin{align*}
    \Bar   & : \text{augmented operad} \rightarrow \text{cooperad}               \\
    \Cobar & : \text{coaugmented cooperads} \rightarrow \text{augmented operads} \\
  \end{align*}

  These otherwise noted in literature as the bar and cobar \emph{constructions}, which can be explained by the fact that $\Bar\P$ and $\Cobar\C$ are chain complexes, that induce differentials $d_\Bar$ and $d_\Cobar$ respectively.

  In terms of the tree functor, $\Bar$ and $\Cobar$ can concretely be defined as

  \begin{align*}
    \Bar\P   & := \T^c(s\overline{\P}, d = d_\Bar + d_\P)      \\
    \Cobar\P & := \T(s^{-1}\overline{\C}, d = d_\Cobar + d_\C) \\
  \end{align*}

  where,

  \begin{align*}
    \overline{\P} & := \ker(\P \rightarrow \I)   & \text{termed augmentation ideal}     \\
    \overline{\C} & := \coker(\I \rightarrow \C) & \text{termed coaugmentation coideal} \\
  \end{align*}

  $d_\Bar$ and $d_\Cobar$ remain to be explained. Let $\P = (\P, \gamma, \eta, \epsilon)$ be an augmented operad. The \Sbb-module $\P$ is naturally isomorphic to $\P = \I \oplus \overline{\P}$. The bar construction $\Bar\P$ is a dg cooperad, whose underlying space is the cofree cooperad $\T^c(s\overline{\P})$ on the suspension of $\overline{\P}$. Since $\T^c(s\overline{\P})$ is a cofree cooperad, there exists a unique coderivation extending the composition on $\P$, which induces the differential $d_\Bar$, $\T^c(s\overline{\P}) \rightarrow \T^c(s\overline{\P})$. Similarly, the differential $d_\Cobar$ is induced by $\T(s^{-1}\overline{\C}) \rightarrow \T(s^{-1}\overline{\C})$. Here, $s$ and $s^{-1}$ refer to \emph{degree suspension} and \emph{degree desuspension}; refer to \cite{Loday12} Section 6.5 for more.
\end{definition}

\begin{theorem}[Relationship between $\Bar$ and \Cobar]
  The functors $\Bar$ and $\Cobar$ are related by the following adjunction.

  \begin{equation*}
    \begin{tikzcd}
      \Cobar \arrow[r, bend left] \arrow[r, phantom, "\bot" description] & \Bar \arrow[l, bend left]
    \end{tikzcd}
  \end{equation*}
\end{theorem}

Twisting morphisms are exactly what induce maps to bar constructions, and maps from cobar constructions.

\begin{theorem}[$\Bar-\Cobar$ resolution]
  The unit and counit of the above adjunction

  \begin{align*}
    \Cobar\Bar\P \rightarrow \P \\
    \C \rightarrow \Bar\Cobar\C \\
  \end{align*}

  are quasi-isomorphisms of operads and cooperads, respectively.
\end{theorem}

\begin{theorem}[Fundamental theorem of \Tw\label{thm:tw}]
  Given a $\alpha \in \Tw(\C, \P)$, the following assertions are equivalent.

  \begin{enumerate}
    \item[(a)] $\P \circ_\alpha \C$ is acyclic.
    \item[(b)] $\C \circ_\alpha \P$ is acyclic.
    \item[(c)] $\C \rightarrow \Bar\P$ is a quasi-isomorphism.
    \item[(d)] $\Cobar\C \rightarrow \P$ is a quasi-isomorphism.
  \end{enumerate}

  In this case, the morphism is termed \emph{Koszul}.
\end{theorem}

\subsection{Koszul duality}
In this section, we restrict our attention to studying the common case of \emph{quadratic operads}, since it turns out that the operads of interest, $\Ass$, $\Com$, and $\Lie$, corresponding to associative, commutative, and Lie algebras, do indeed turn out to be quadratic. In fact, $\Ass$, $\Com$, and $\Lie$ are \emph{binary quadratic operads}, which is to say that they are generated by binary operations and quadratic relations.

\begin{definition}[Quadratic data]
  Let us first recall the \emph{minimal model} of an operad \P, as defined in \ref{def:minmodel-op}: a given $\M = \T(E)$ is \emph{minimal} when $d(E) \subset \T(E)^{\geq 2}$. To get the definition of a quadratic operad, we restrict $d(E)$ further to $d(E) \subset \T(E)^{(2)}$. Here, $\T(E)^{(2)}$ is the sub-\Sbb-module of the free operad \T, which is spanned by two elements of $E$.

  A quadratic operad or cooperad can be constructed with two pieces of data, say $(E, R)$, where elements of $E$ are termed \emph{generating operations} and elements of $R$ are termed \emph{relators}. The quadraticity hypothesis says that the relators are made up of elements involving exactly two compositions. We will, at the end of the section, see the example of $\Ass$ and \Com, which are \emph{binary quadratic operads}; the extra ``binary'' qualifier says that there the generating operations consist entirely of binary operations, in which case, all relations involve three variables.
\end{definition}

\begin{definition}[Quadratic operad]
  A quadratic operad $\P$ is the free operad over $E$ quotiented out by the ideal generated by $R$.

  \begin{equation*}
    \P(E, R) := \T(E)/(R)
  \end{equation*}

  such that the composite

  \begin{equation*}
    \begin{tikzcd}
      R \arrow[r, rightarrowtail] & \T(E) \arrow[r, twoheadrightarrow] & \P
    \end{tikzcd}
  \end{equation*}

  is zero, and the following diagram commutes

  \begin{equation*}
    \begin{tikzcd}
      R \arrow[r, rightarrowtail] & \T(E) \arrow[r, twoheadrightarrow] \arrow[rd, rightarrowtail] & \P \arrow[d, dashed, "\exists!" description] \\
      & & \P(E, R)
    \end{tikzcd}
  \end{equation*}
\end{definition}

\begin{example}[Explication of quadratic operad]
  Let $E = \{0, E(1), E(2), E(3), E(4), \ldots\}$. Then, $\T(E) = \oplus_k \T(E)^{(k)}$, where

  \begin{equation*}
    \begin{matrix}
      \T(E)^{(0)} & = & (0, \K\id, 0, 0, \ldots)                                          & = & \I \\
      \T(E)^{(1)} & = & (0, E(1), E(2), E(3), \ldots)                                     & = & E  \\
      \T(E)^{(2)} & = & (0, E(1)^{\otimes 2}, E(2)^{\otimes 2}, E(3)^{\otimes 2}, \ldots) &   &    \\
    \end{matrix}
  \end{equation*}

  The quotient $\P(E, R)$ of $\T(E)$ is such that

  \begin{align*}
    \P(E, R)^{(0)} & = \I                                 \\
    \P(E, R)^{(1)} & = E                                  \\
    \P(E, R)^{(2)} & = (0, E(1)^{\otimes 2}/R(1), \ldots) \\
  \end{align*}
\end{example}

\begin{definition}[Quadratic cooperad]
  A quadratic cooperad $\C$ is the sub-cooperad of the co-free cooperad $\T^c(E)$ such that the composite

  \begin{equation*}
    \begin{tikzcd}
      \C \arrow[r, rightarrowtail] & \T^c(E) \arrow[r, twoheadrightarrow] & \T^c(E)^{(2)}/R
    \end{tikzcd}
  \end{equation*}

  is zero, and the following diagram commutes

  \begin{equation*}
    \begin{tikzcd}
      \C \arrow[r, rightarrowtail] \arrow[d, dashed, "\exists!" description] & \T^c(E) \arrow[r, twoheadrightarrow] & \T^c(E)^{(2)}/R \\
      \C(E, R) \arrow[ru, rightarrowtail] & &
    \end{tikzcd}
  \end{equation*}
\end{definition}

\begin{example}[Explication of quadratic cooperad]
  Contrasting with the explication of a quadratic operad, we see that $\T^c(E)$ has the same underlying \Sbb-module as $\T(E)$, and $\C(E, R)$ is such that

  \begin{align*}
    \C(E, R)^{(0)} & = \I                            \\
    \C(E, R)^{(1)} & = E                             \\
    \C(E, R)^{(2)} & = (0, R(1), R(2), R(3), \ldots) \\
  \end{align*}
\end{example}

We now take a second look at many of notions introduced in the previous section, in the special case of quadratic operads.

\begin{definition}[\Pza]
  We define the \emph{Koszul dual cooperad} of a quadratic operad, \Pza~\footnote{prononouced "anti-shriek"}, as

  \begin{equation*}
    \Pza(E, R) := \C(E[1], R[2])
  \end{equation*}

  By definition,

  \begin{equation*}
    (\Pza)^{\text{!`}} = \P
  \end{equation*}
\end{definition}

\begin{definition}[\Pz]
  The \emph{Koszul dual operad}, \Pz~\footnote{pronounced "shriek"}, of a quadratic operad $\P$ is the arity-wise linear dual of \Pza, where the square brackets indicate degree shift:

  \begin{equation*}
    \Pz(n) := \Hom_\K(\Pza(n), \K)[1 - n]
  \end{equation*}
\end{definition}

\begin{lemma}[Self-dualizing property of shriek]
  For quadratic operad \P, generated by a reduced \Sbb-module of finite dimension in each arity,

  \begin{equation*}
    (\Pz)^{!} = \P
  \end{equation*}

  In particular, the finite dimensionality condition ensures that $(E^\vee)^\vee = E$.
\end{lemma}

\begin{example}[Relationships between \Ass, \Com, and \Lie]
  In summary,

  \begin{align*}
    \Ass^{!} & = \Ass \\
    \Com^{!} & = \Lie \\
    \Lie^{!} & = \Com \\
  \end{align*}
\end{example}

\begin{remark}[$\Bar$ and $\Cobar$ for quadratic operad]
  Since the internal differential of $E$ is trivial,

  \begin{enumerate}
    \item[(\P)] The inclusion $\Pza \rightarrowtail \Bar\P$ induces an isomorphism $\Pza \cong H^0(\Bar^\bullet\P)$.
    \item[(\C)] The surjection $\Cobar\C \twoheadrightarrow \Cza$ induces an isomorphism $H_0(\Cobar_\bullet\C) \cong \C$.
  \end{enumerate}
\end{remark}

\begin{remark}[Natural twisting morphism]
  Given quadratic data $(E, R)$, since $\P(E, R)^{(1)} = E$ and $\C(E, R)^{(1)} = E$, we can define a morphism $\kappa$ as the composite

  \begin{equation*}
    \C(E[1], R[2]) \twoheadrightarrow sE \rightarrow E \rightarrowtail \P(E, R)
  \end{equation*}

  Notice that $\kappa \star \kappa = 0$, and hence $\kappa$ is a twisting morphism.
\end{remark}

\begin{definition}[Koszul complex]
  One complex associated to $\kappa$,

  \begin{equation*}
    \Pza \circ_\kappa \P := (\Pza \circ \P, d_\kappa)
  \end{equation*}

  is termed as the \emph{(right) Koszul complex} associated to quadratic operad \P. We can define another completely different complex corresponding to $\P \circ_\kappa \Pza$, and this would be called the \emph{left Koszul complex}.
\end{definition}

\begin{terminology}[Koszuality of an operad]
  A quadratic operad $\P$ is termed \emph{Koszul} if and only if its associated Koszul complex is acyclic.
\end{terminology}

\begin{lemma}[Minimal model of Koszul operad]
  The minimal model of a Koszul operad $\P$ is $\Cobar\Pza$.
\end{lemma}

\begin{lemma}[Koszuality and duality]
  Finitely-generated operad $\P$ is Koszul if and only if $\Pz$ is Koszul.
\end{lemma}

\begin{terminology}[\Pinf-algebra]
  Given Koszul operad \P, by definition, a \emph{homotopy \P-algebra}, or a \Pinf-algebra, is an algebra over \Pinf, which is defined as

  \begin{equation*}
    \Pinf := \Cobar\Pza
  \end{equation*}
\end{terminology}

\begin{discussion}[Canonical cofibrant resolution]
  Recalling \ref{thm:tw}, we infer that for a Koszul operad \P, $\Pza \rightarrow \Bar\P$ and $\Cobar\Pza \rightarrow \P$ are quasi-isomorphisms, and they can be thought of as resolutions of \P.

  Indeed, with a bit of background on the model structure on \P~\footnote{Otherwise known as the Hinich model structure \cite{Hinich97}.}, which we refrain from discussing here, it turns out that $\Pinf$ is the \emph{canonical cofibrant resolution} of \P.
\end{discussion}

\begin{example}[Duals of $\Com$ and \Lie]
  As mentioned earlier, $\Com$ and $\Lie$ are binary quadratic operads. We take this opportunity to work out their duals.

  The quadratic operad $\Com = \P(E, R)$ is generated by binary symmetric operations in $\mu$, i.e. $E = \K\mu$. Let us denote by $\mu_\mathsf{I}$, $\mu_\mathsf{II}$ and $\mu_\mathsf{III}$ the elements $\mu \circ_\mathsf{I} \mu$, $\mu \circ_\mathsf{II} \mu$, and $\mu \circ_\mathsf{III} \mu$ in $\T(\K\mu)(3)$. The space of relations $R$ admits the following \K-linear basis made up of two elements, $\mu_\mathsf{II} - \mu_\mathsf{I}$ and $\mu_\mathsf{III} - \mu_\mathsf{II}$.

  The Koszul dual is generated by the space $E^\vee \cong \K\nu$, where $\nu^{(12)} = -\nu$. Again, let us denote by $\nu_\mathsf{I}$, $\nu_\mathsf{II}$, and $\nu_\mathsf{III}$, the elements of $\nu \circ_\mathsf{I} \nu$, $\nu \circ_\mathsf{II} \nu$, and $\nu \circ_\mathsf{III} \nu$ of $\T(\K\nu)(3)$, respectively. The orthogonal space of $R$ under scalar product $\langle -, - \rangle$ has dimension with basis $\nu_\mathsf{I} + \nu_\mathsf{II} + \nu_\mathsf{III}$. Denoting skew-symmetric generating operation $\nu$ by bracket $[-, -]$, the latter is nothing but the Jacobi identity $[[a, b], c] + [[b, c], a] + [[c,a], b] = 0$. Hence, $\Com^{!} = \Lie$.
\end{example}

\begin{summary}
  Since all our operads of interest are Koszul, we can safely assume that Koszul duality works, and indeed that, a homotopy operadic algebra is well-defined. In the final section, we will shift our attention exclusively to \Pinf-algebras.
\end{summary}

\subsection{\texorpdfstring{\Pinf-algebras}{Homotopy operadic algebras}}
All that we've studied so far culminate into this final section containing the seminal results of operadic theory.

The main topic of study is a theorem central to study of homotopy \P-algebras, the \emph{Rosetta Stone}, which provides isomorphisms between four different definitions of a \Pinf-algebra. We state it as-is, with some commentary.

\begin{theorem}[Rosetta Stone]
  The Rosetta Stone establishes that four different explications of the \Pinf- structure on module $A$ are equivalent.

  \begin{equation*}
    \Hom(\Cobar\Pza, \End{A}) \cong \Tw(\Pza, \End{A}) \cong \Hom(\Pza, \Bar\End{A}) \cong \mathsf{Codiff}(\Pza(A))
  \end{equation*}
\end{theorem}

\begin{discussion}
  \begin{equation*}
    \Cobar\Pza \tilde{\longrightarrow} \P \longrightarrow \End{A}
  \end{equation*}

  can be interpreted in terms of the $\Tw$. A morphism of operads from $\Cobar\Pza$ to $\End{A}$ is equivalent to the twisting morphism in the convolution algebra

  \begin{equation*}
    \g := \Hom_\Sbb(\Pza, \End{A})
  \end{equation*}

  A homotopy \P-algebra structure on module $A$ is equivalent to the twisting morphism $\Tw(\Pza, \End{A})$.
\end{discussion}

\begin{lemma}[Characterization of \P-algebras]
  \P-algebras are characterized among \Pinf-algebras as particular solutions to the Maurer-Cartan equation. In particular, a \Pinf-algebra is a \P-algebra if its twisting morphism is concentrated in weight (number of vertices in the tree) $1$.
\end{lemma}

\begin{definition}[\8-isomorphism and \8-quasi-isomorphism]
  An \8-morphism of $\Pinf$ algebras is given by:
  \begin{equation*}
    f: (\Pza(A), d_\sigma) \rightarrow (\Pza(B), d_\psi)
  \end{equation*}

  It is an \8-isomorphism if its first component $f_{(0)}: A \rightarrow B$ is an isomorphism.

  An \8-quasi-isomorphism $f$ is called so if its first component $f_{(0)}: A \rightarrow B$ is a quasi-isomorphism.
\end{definition}

\begin{theorem}[Homotopy Transfer Theorem]
  When working over a field, the homotopy $H(A)$ can be made a deformation retract of $A$. It enables us to transfer the $\Pinf$-algebra structure from $A$ to $H(A)$.

  Let $(V, d_V)$ be the homotopy retract of the data of $(W, d_W)$:

  \begin{equation*}
    \begin{tikzcd}
      \arrow[loop left, "{h}"](W, d_W) \arrow[r, shift left=1, "{p}"] & \arrow[l, shift left=1, "{i}"] (V, d_V)
    \end{tikzcd}
  \end{equation*}
  Let $\P$ be a Koszul operad, and let $(V, d_V)$ be the homotopy retract of $(W, d_W)$. Any $\Pinf$-algebra structure on $W$ can be transferred to a $\Pinf$-algebra structure on $V$ such that $i$ extends to an \8-quasi-isomophism.

  To prove it, we use the third definition of a $\Pinf$-algebra supplied in the Rosetta Stone.

  \begin{equation*}
    \Hom(\Cobar\Pza, \End{A}) \cong \Tw(\Pza, \End{A}) \cong \Hom(\Pza, \Bar\End{A})
  \end{equation*}
\end{theorem}

\section{Derived deformation theory\label{sec:derived-deformation}}
In the first two sections, we introduce and motivate derived deformation theory. The final section introduces the Lurie-Pridham equivalence, and can be read as an introduction to derived deformation theory, as it appears in \cite{Lurie18}.

\subsection{Classical deformation theory of associative algebras}
In order to understand precisely what a ``deformation'' is, we provide some intuitions from its historical roots: the classical deformation theory of associative algebras.

\begin{definition}[Hochschild cohomology]
  The Hochschild cohomology of an associative algebra $A$ with coefficients in module $M$ is given by:

  \begin{equation*}
    \begin{tikzcd}
      0 \arrow[r, "{\delta_\Hoch}"] & M \arrow[r, "{\delta_\Hoch}"] & C^1_\Hoch(A, M) \arrow[r, "{\delta_\Hoch}"] & \ldots \arrow[r, "{\delta_\Hoch}"] & C^n_\Hoch(A, M) \arrow[r, "{\delta_\Hoch}"] & \ldots
    \end{tikzcd}
  \end{equation*}

  where $C^n_\Hoch(A, M)$ is the space of $n$-multilinear maps from $A$ to $M$.

  \begin{equation*}
    C^n_\Hoch(A, M) := \Hom_\K(A^{\otimes n}, M)
  \end{equation*}
\end{definition}

\begin{definition}[Formal deformation]
  A formal deformation of associative algebra $A$ is defined as a deformation along the formal power series $\K\llbracket t \rrbracket$, given by the family

  \begin{equation*}
    \{\mu_i : A \otimes A \rightarrow A \mid i \in \mathbb{N}\}
  \end{equation*}

  where each $\mu_i$ corresponds to the $t^i$ deformation, satisfying the following:

  \begin{enumerate}
    \item[($D_0$)] $\mu_0(a, b) = ab$.
    \item[($D_1$)]
      \begin{align*}
        \mu_1(a, \mu_0(b, c)) & + \mu_0(a, \mu_1(b, c)) = \\
        \mu_1(\mu_0(a, b), c) & + \mu_0(\mu_1(a, b), c)
      \end{align*}
    \item[($D_k$)] $\Sigma_{i + j = k} \mu_i(a, \mu_j(b, c)) = \Sigma_{i + j = k} \mu_i(\mu_j(a, b), c)$
  \end{enumerate}

  where $a, b, c \in A$.
\end{definition}

\begin{remark}[An observation on $D_1$]
  Observe that $D_1$ reads as

  \begin{equation*}
    a\mu_1(b, c) - \mu_1(ab, c) + \mu_1(a, bc) - \mu_1(a, b)c = 0
  \end{equation*}

  This says that $\mu_1 \in \Hom_\K(A^{\otimes 2}, A)$ is a Hochschild cocycle. In other words, $\delta_\Hoch(\mu_1) = 0$. This motivates our next definition.
\end{remark}

\begin{definition}[Infinitesimal deformation]
  An infinitesimal deformation $\epsilon$, concentrated only at degree $0$, of associative algebra $A$ over field $\K$ is a $D$-deformation of $A$, where

  \begin{equation*}
    D := \K[\epsilon]/\epsilon^2
  \end{equation*}

  is a \emph{local Artin ring} [\ref{def:art-rng}] of dual numbers.
\end{definition}

\begin{theorem}[Significance of Hochschild cochain]
  There is a one-to-one correspondence between equivalence classes of infinitesimal deformations of $A$ and Hochschild cohomology $H_\Hoch^2(A, A)$.
\end{theorem}

\begin{intuition}[Maurer-Cartan equation in associative algebras]
  Let $A$ be an associative algebra and multiplication operator $\mu: A \otimes A \rightarrow A$ and $\nu \in C_\Hoch(A, A)$ be a Hochschild 2-cochain. Then, $\mu + \nu$ is associative if and only if

  \begin{equation*}
    \delta_\Hoch(\nu) + \frac{1}{2}[\nu, \nu] = 0
  \end{equation*}

  This is indeed an incarnation of the Maurer-Cartan equation, and this gives us the first hint that there is a Lie algebra hiding in every deformation problem.
\end{intuition}

\subsection{Motivation from classical geometric deformation theory}
We now present the journey from (classical) deformation theory in algebraic geometry to derived deformation theory, in derived algebraic geometry.

In algebraic geometry, the deformation functor

\begin{equation*}
  \mathsf{Def}_{X_0}: \Art_\K \rightarrow \Set
\end{equation*}

sends an augmented Artin algebra $\A$ to isomorphism classes $(X, u)$, where $X$ is a scheme and $u: X \otimes_\A \K \cong X_0$ is an isomorphism of \emph{algebraic varieties}.

It was then realized that this deformation problem was controlled by a dg-Lie algebra $L$, and we form the set $\MC(L \otimes_\K m)$, where $m$ is the maximal ideal of an augmented Artin algebra $\A$.

\begin{equation*}
  X_L: \Art_\K \rightarrow \Set
\end{equation*}

The ring of functions formed by $X_L$ can be identified with $H^0(\hat{C}^*(L))$ where $\hat{C}^*: \op{\mathsf{dgLie}}_\K \rightarrow \mathsf{cdga}_\K$ is termed \emph{Chevalley complex}.

In the above picture, the \emph{tangent space} is a simple $H^1(L)$, and $H^2(L)$ gives us one obstruction to the formal smoothness. $H^i(L)$ for $i \leq 0$ gives us increasing orders of symmetry, and the information contained in $H^i(L)$ for $i \geq 2$ is most interesting, and marks the beginning of derived deformation theory; indeed, it gives higher order defects, and there is no way to extract this information from the classical picture.

We proceed to refine the above picture.

\begin{equation*}
  \F: \Art_\K \rightarrow \Set
\end{equation*}

Functor $\F$ sends (augmented) Artin algebras to \emph{equivalence classes} of moduli spaces. The tangent space of $\F$, formerly given by $\F(\K[\epsilon]/\epsilon^2)$, can now be thought of as a \emph{tangent vector}, $\F(\K[x]/x^3) \rightarrow \F(\K[x]/x^2)$. It then determines the following \emph{homotopy cartesian} square

\begin{equation*}
  \begin{tikzcd}
    \K[x]/x^3 \arrow[r] \arrow[d] & \K \arrow[d] \\
    \K[x]/x^2 \arrow[r] & \K[\epsilon_{-1}]/{\epsilon_{-1}}^2 \\
  \end{tikzcd}
\end{equation*}

where $\F(\K[\epsilon_{-1}]/{\epsilon_{-1}}^2)$ corresponds to the classical deformation problem, and $\epsilon_{-1}$ is in cohomological degree $-1$.

If $\F$ comes from the restriction of a functor defined on the dg-Artin algebra, we can study the smoothness defects of $\F$ at $x$ using the above homotopy cartesian square.

The \emph{equivalence classes} mentioned in the above picture is also vague; using machinery we will detail in \ref{sec:goldman}, the equivalence classes are then identified as weak equivalences of simplicial sets. Further, the functor $\F$ can be identified as the \emph{formal spectrum} of augmented Artin algebra \A, denoted \Spf(\A).

\begin{equation*}
  \Spf(\A): \mathsf{dgArt}_\K \rightarrow \SSet
\end{equation*}

Having to work with weak equivalences of simplicial sets suggests that the picture can be lifted to the level of abstraction of \8-categories, which we describe in the next section.

\subsection{The Lurie-Pridham equivalence\label{sec:lurie-pridham-equiv}}
The Lurie-Pridham equivalence, in its elaborated form, involves a lot of machinery. We supply the formal definitions as-is, and supply intuitions.

\begin{convention}
  Lurie's work is done in the setting of \8-operads, and in particular, all algebras mentioned are \8-operadic algebras. However, we will not concern ourselves with the details of this.
\end{convention}

\begin{definition}[Deformation context]
  A deformation context is a pair $(\A, \{E_\alpha\}_{\alpha \in T})$, where $\A$ is a \emph{presentable \8-category} [\ref{int:pres-inf}], and $\{E_\alpha\}_{\alpha \in T}$ is a set of \emph{spectrum objects} of \A, denoted $\Sp(\A)$ [\ref{def:Sp}].
\end{definition}

\begin{intuition}[$\{E_\alpha\}_{\alpha \in T}$]
  This can be thought of as $\K[\epsilon_i]/({\epsilon_i}^2)$ in the classical setting, where $\epsilon$ is concentrated in non-negative degrees; each $\epsilon_i$ is a shift by $\mathsf{[i]}$, corresponding to the loop functor $\Omega$ in stable homotopy theory.
\end{intuition}

\begin{definition}[Formal moduli problem]
  Let $(\A, \{E_\alpha\}_{\alpha \in T})$ be a deformation context. A formal moduli problem is defined as a functor from an \emph{Artin algebra} [\ref{def:art-alg}] to the \emph{\8-category of spaces} [\ref{def:inf-space}].
  \begin{equation*}
    \A^\Art \rightarrow \S
  \end{equation*}

  satisfying:

  \begin{enumerate}
    \item[(a)] Space $X(*)$ is contractible, where $*$ denotes final object of \A.
    \item[(b)] Let $\sigma$:

      \begin{equation*}
        \begin{tikzcd}
          A' \arrow[r]\arrow[d] & B' \arrow[d, "{\phi}"] \\
          A \arrow[r] & B \\
        \end{tikzcd}
      \end{equation*}

      be a diagram in $\A^\Art$. If $\sigma$ is a pullback diagram and $\phi$ is small, then $X(\sigma)$ is a pullback diagram in \S.
  \end{enumerate}
\end{definition}

\begin{definition}[Deformation problem\label{def:defaxioms}]
  A deformation problem is a functor $\Dfrak : \op{\A} \rightarrow \B$, where $\A$ and $\B$ are some algebras, satisfying the following list of axioms.

  \begin{enumerate}
    \item[(D1)] $\B$ is a \emph{presentable \8-category} [\ref{int:pres-inf}].
    \item[(D2)] $\Dfrak$ admits a left adjoint $\Dfrak' : \B \rightarrow \op{\A}$.
    \item[(D3)] There exists a full subcategory $\B_0 \subseteq \B$ satisfying:

      \begin{enumerate}
        \item[(a)] For every $K \in B_0$ the unit map $K \rightarrow \Dfrak\Dfrak'K$ is an equivalence.
        \item[(b)] $\B_0$ contains the initial object $\0 \in \B$. Then, it follows from (a) that $\0 \simeq \Dfrak\Dfrak'\0 \simeq \Dfrak(*)$ where $*$ is the initial object of $\A$.
        \item[(c)] For every $\alpha \in T$, and every $n \geq 1$, there exists an object $K_{\alpha, n} \in B_0$, and an equivalence $\Omega^{\8 - n} E_\alpha \simeq \Dfrak'K_{\alpha, n}$.
        \item[(d)] For every pushout:

          \begin{equation*}
            \begin{tikzcd}
              K_{\alpha, n} \arrow[r] \arrow[d] & K \arrow[d] \\
              \0 \arrow[r] & K' \\
            \end{tikzcd}
          \end{equation*}

          where $\alpha \in T$, $n \in \mathbb{N}$, if $K \in B_0$, then $K' \in B_0$.
      \end{enumerate}

    \item[(D4)] For every $\alpha \in T$, $e_\alpha : \B \rightarrow \Sp(\S)$ preserves \emph{small sifted colimits} [\ref{int:siftlim}]. Moreover, a morphism $f \in \B$ is an equivalence if and only if each $e_\alpha(f)$ is an equivalence of spectra.
  \end{enumerate}
\end{definition}

\begin{intuition}[Interpretation of the axioms\label{int:defaxioms}]
  Let us try to approximate what (D1)-(D4) mean.

  \begin{enumerate}
    \item[(1)] $\Dfrak$ and $\Dfrak'$ form an adjoint pair. They're almost equivalences, but for two caveats:
      \begin{enumerate}
        \item[(a)] The counit map is not an equivalence.
        \item[(b)] The unit map is an equivalence only for some ``good'' subcategory $\B_0 \subseteq \B$.
      \end{enumerate}
    \item[(2)] $\Dfrak'$ is homotopy inverse to $\Dfrak$ on some ``good'' subcategory $\B_0 \subseteq \B$. In particular, (D3d) is dual to the condition on Artin algebras [\ref{def:artin-palg}].
    \item[(D4)] $e_\alpha$ admits a left adjoint, and $\B \simeq \mathsf{LMod}_U(Sp^T)$ is an equivalence of \8-categories. In other words, we can think of objects of $\B$ being determined by a collection of spectra indexed by $T$ equipped with the additional structure of \emph{left action by monad U} [\ref{int:barrbeck}].
  \end{enumerate}

  We will return to this once again, in the \ref{sec:op-deformation}.
\end{intuition}

\begin{theorem}[Lurie-Pridham]
  Given a deformation problem $\Dfrak: \op{\A} \rightarrow \B$, it is related to the class to the class of formal moduli problems \FMP, as follows.

  \begin{equation*}
    \FMP^\A \simeq \B
  \end{equation*}

  is an equivalence of \8-categories.
\end{theorem}

\subsection{The tangent complex}
The tangent complex is an enhancement of the tangent vector space into a chain complex, such that it is defined at every point in the space and globally as a bundle. For smooth points, tangent complexes are quasi-isomorphic to tangent spaces. The most important property of tangent complexes is that they have a natural Lie algebra structure; at smooth points, it is Abelian. The whole formal neighborhood of a singular point can be reconstructured from this structure.

Let $X$ be an algebraic variety over the field $\mathbf{C}$ of complex numbers, and let $x: X \rightarrow \Spec(\mathbf{C})$ be a point of $X$. Then a tangent vector at point $x$ is the dotted arrow rendering the following diagram commutative.

\begin{equation*}
  \begin{tikzcd}
    \Spec(\mathbf{C}) \arrow[d] \arrow[r, "{x}"] & X \arrow[d] \\
    \Spec(\mathbf{C}[\epsilon]/\epsilon^2) \arrow[r] \arrow[ur, dashed] & \Spec(\mathbf{C})
  \end{tikzcd}
\end{equation*}

The collection of tangent vectors at $x$ comprise a vector space $\Tbb_x$. Our goal is to generalize the construction of $\Tbb_x$ in the setting of an arbitrary formal moduli problem.

\begin{definition}[Tangent space]
  Let $(\A, \{E_\alpha\}_{\alpha \in T})$ be a deformation context and let $Y: \A^\Art \rightarrow \S$ be a formal moduli problem. For each $\alpha \in T$, the tangent space of $Y$ at $\alpha$ is the space $Y(\Omega^\8 E_\alpha)$.
\end{definition}

\begin{lemma}[Characterization of spectrum objects]
  If $\Cscr$ is an \8-category which admits finite limits, the spectrum objects of \Cscr, or $\Sp(\Cscr)$, is defined to be the full subcategory of $\Hom(\Sfin, \Cscr)$, spanned by those functors which are reduced and excisive. [\ref{not:Sfin}]
\end{lemma}

Let $(\A, \{E_\alpha\}_{\alpha \in T})$ be a deformation context. Then, we identify, for each $\alpha \in T$, $E_\alpha \in \Sp(\Cscr)$ with the corresponding functor $\Sfin \rightarrow \Cscr$. It can be shown that, for pointed space $K$, $E_\alpha(K)$ is Artinian. It can further be shown that the composite functor

\begin{equation*}
  \Sfin \overset{E_\alpha}{\rightarrow} \A^\Art \overset{Y}{\rightarrow} \S
\end{equation*}

is reduced and excisive.

These observations motivate the following definition.

\begin{definition}[Tangent complex]
  Let $(\A, \{E_\alpha\}_{\alpha \in T})$ be a deformation context, and $Y: \A^\Art \rightarrow \S$ be a formal moduli problem. Then, the composite functor

  \begin{equation*}
    \Sfin \overset{E_\alpha}{\rightarrow} \A^\Art \overset{Y}{\rightarrow} \S
  \end{equation*}

  denoted by $Y(E_\alpha)$, is an object of \8-category $\Sp(\S)$. We will call $Y(E_\alpha)$ the \emph{tangent complex} of $Y$ at $\alpha$.
\end{definition}

\begin{remark}[Tangent spectrum]
  The terminology introduced in the above definition is potentially misleading, as $Y(E_\alpha)$ is a spectrum. However, for instance, when $(\A, \{E_\alpha\}_{\alpha \in T}) = (\mathsf{CAlg}^{\mathsf{aug}}_\K, \{E\})$, for field $\K$, the tangent spectrum $Y(E)$ admits the structure of a \K-module spectrum, and can hence be identified with chain complexes of vector spaces over \K. See also: [\ref{thm:dold-kan}]
\end{remark}

\section{Operadic deformation theory\label{sec:op-deformation}}
This part is dedicated to re-interpreting \ref{sec:derived-deformation} in the setting of operads, for which the reader should be fully prepared, after having digested \ref{sec:operads}. All original results will be marked with a $\dagger$.

\subsection{The Goldman-Millson theorem\label{sec:goldman}}
The aim of this interlude is to understand the mechanics of \emph{simplicial enrichment} of a Lie algebra \g, to better equip ourselves to understand the next sections.

We begin with the well-known concept of the Maurer-Cartan set of elements of Lie algebra \g, \MC(\g) [\ref{def:mc-elem}]. It has the defining property that, a field $\g^1 \rightarrow \g^1$ with a \emph{flow}, preserves the Maurer-Cartan set of elements of \g. We say that $\MC(\g)$ is \emph{flow-invariant}, and we obtain \emph{gauge equivalences} in \g, which we denote as $\sim_\mathsf{gauge}$. We can now define

\begin{definition}[Maurer-Cartan moduli space]
  The \emph{moduli space} of Maurer-Cartan elements is defined as

  \begin{equation*}
    \overline{\MC}(\g) := \MC(\g)/\sim_\mathsf{gauge}
  \end{equation*}

  In fact,

  \begin{equation*}
    \overline{\MC}: \{\g \in \Lie \mid \text{$\g$ is complete with respect to its canonical filtration [\ref{term:complete-canonical}]}\} \rightarrow \Set
  \end{equation*}
\end{definition}

\begin{theorem}[Goldman-Millson]
  Let \g, $\h$ be \emph{complete Lie algebras} [\ref{def:complete-lie}], and $\phi: \g \rightarrow \h$ be a \emph{filtered quasi-isomorphism} [\ref{def:filtered-qiso}]. Then,

  \begin{equation*}
    \overline{\MC}(\phi): \overline{\MC}(\g) \rightarrow \overline{\MC}(\h)
  \end{equation*}

  is a bijection. See \cite{Nicoud18} for more on this.
\end{theorem}

The Goldman-Millson theorem is but an avatar of a ``higher'' result. This result encodes higher equivalences into \MC, not just gauge equivalences.

Consider a nilpotent Lie algebra \g, and tensor it with the \emph{Sullivan algebra} $\Omega_\bullet$ [\ref{def:sullivan-alg}]. Notice that $\g \otimes \Omega_n$ is again a Lie algebra by

\begin{equation*}
  [x \otimes \alpha, y \otimes \beta] = [x, y] \otimes \alpha\beta
\end{equation*}

\begin{definition}[Simplicial enrichment of Lie algebra\label{def:simpenrich-lie}]
  Given a complete Lie algebra $(\g, F_\bullet\g)$, we can \emph{simplicially enrich} it by tensoring it with the Sullivan commutative algebra $\Omega_\bullet$, as follows:

  \begin{equation*}
    \g \otimes \Omega_\bullet := \lim_n(\g/F_n\g \otimes \Omega_\bullet)
  \end{equation*}
\end{definition}

\begin{definition}[Hinich-Getzler \8-groupoid]
  Let $\g$ be a complete Lie algebra. The \emph{Hinich-Getzler \8-groupoid} is defined as the simplicial set

  \begin{align*}
    \MC_\bullet     & : \textsf{dg-Lie} \rightarrow \SSet \\
    \MC_\bullet(\g) & := \MC(\g \otimes \Omega_\bullet)   \\
  \end{align*}

  As hinted earlier, indeed

  \begin{equation*}
    \pi_0 \MC_\bullet(\g) = \overline{\MC}(\g)
  \end{equation*}
\end{definition}

\begin{theorem}[Fibrancy of $\MC_\bullet$\label{thm:mc-infcat}]
  The functor $\MC_\bullet$ has image in the full subcategory of \Kan. See \cite{Hinich96} for more.
\end{theorem}

Finally, we state the more modern version of the Goldman-Millson theorem.

\begin{theorem}[Dolgushev-Rogers]
  Let \g, $\h$ be complete Lie algebras, and $\phi: \g \rightarrow \h$ be a filtered quasi-isomorphism. Then,

  \begin{equation*}
    \MC_\bullet(\phi): \MC_\bullet(\g) \rightarrow \MC_\bullet(\h)
  \end{equation*}

  is a homotopy equivalence of simplicial sets. See \cite{Dolgushev15} for more.
\end{theorem}

\begin{summary}
  The functor $\MC$ has several concrete definitions, and its image is an \8-category, in its most general form. To avoid encumbering ourselves with notation, we will henceforth refer to all variations as \MC, letting the context resolve the ambiguity.
\end{summary}

\subsection{Main results}
\begin{convention}
  The setting is homotopy invariant, and in particular, \P-algebras and  \Pinf-algebras are indistinguishable in this setting.
\end{convention}

\begin{convention}
  All (co)operads are assumed to be (co)augmented, and (filtered) cofibrant as left \K-modules.
\end{convention}

\begin{convention}
  The paper \cite{Calaque19} is written in the setting of \emph{colored operads}. In order to avoid encumbering ourselves with notation and fine detail, we will drop the colors, which must be drawn from \K, everywhere.
\end{convention}

\begin{olemma}[\P-algebras form an \8-category]
  The category of \P-algebras, has image in the full subcategory of \Kan.

  \begin{proof}
    In the setting of \8-operads, the machinery of which is detailed in \cite{Lurie12}, this is true by construction. However, we are working in the setting of ordinary operads (with color implicit), there is only a model theoretic structure on it, and there is a simplicial enrichment possible on a ordinary \P-algebra \A just like the case of Lie algebras in field of characteristic $0$ given by \ref{def:simpenrich-lie}, by tensoring with simplicial commutative \emph{Sullivan algebra} $\Omega_\bullet$:

    \begin{equation*}
      \A \otimes \Omega_\bullet
    \end{equation*}
  \end{proof}
\end{olemma}

\begin{definition}[Artin \P-algebra\label{def:artin-palg}]
  The \8-category $\ArtP$ of \emph{Artin \P-algebras} is the smallest full subcategory of the \8-category of \P-algebras such that

  \begin{enumerate}
    \item[(i)] Every trivial algebra $\K[n]$ is Artin, for all $n \in \N$.
    \item[(ii)] Any given Artin \P-algebra \A, the homotopy pullback

      \begin{equation*}
        \begin{tikzcd}
          \A' \arrow[r] \arrow[d] & 0 \arrow[d] \\
          \A \arrow[r] & \K[n] \\
        \end{tikzcd}
      \end{equation*}

      is also Artin, for all $n \geq 1$.
  \end{enumerate}
\end{definition}

\begin{definition}[Formal moduli problem over operad]
  A formal moduli problem is defined as a functor from Artin \P-algebras to the \8-category of spaces.

  \begin{equation*}
    \FMPP: \Art_\P \rightarrow \S
  \end{equation*}

  satisfying the following conditions:

  \begin{enumerate}
    \item[(i)] $F(0) \simeq *$, where $0$ is the zero-algebra.
    \item[(ii)] $F$ sends a pullback diagram in $\ArtP$ of the form

      \begin{equation*}
        \begin{tikzcd}
          \A' \arrow[r]\arrow[d] & 0 \arrow[d] \\
          \A \arrow[r] & \K[n]
        \end{tikzcd}
      \end{equation*}

      to pullback squares of spaces, for each $n \geq 1$.
  \end{enumerate}
\end{definition}

\begin{example}[Formal spectrum of \P-algebra]
  To supply an example of a formal moduli problem, we can associate to every \P-algebra \B, its \emph{formal spectrum} defined as

  \begin{align*}
    \Spf(\B): \ArtP & \rightarrow \S          \\
    \A              & \mapsto \Map_\P(\B, \A) \\
  \end{align*}

  Applying condition (ii) from the definition of a formal moduli problem, for the case $\A = 0$, the formal spectrum can be interpreted as the following sequence of equivalences.

  \begin{equation*}
    \F(\K) \tilde{\longrightarrow} \Omega\F(\K[1]) \tilde{\longrightarrow} \Omega^2\F(\K[2]) \tilde{\longrightarrow} \ldots
  \end{equation*}

  The sequence of spaces forms an $\Omega$-cohomological-spectrum $\Tgt{\F}$.
\end{example}

\begin{olemma}[Homotopy category of formal spectrum]
  The space of morphisms between $\Art_\P$ and $\S$ is defined as the formal spectrum $\Spf$. The homotopy category of this particular formal spectrum (as opposed to the spectrum defined by Lurie, using the machinery of \8-operads) does not possess a triangulated structure [\ref{def:triangulated-cat}].

  \begin{proof}
    $\Spf$ is an $\Omega$-cohomological spectrum, in which the $\pi_*$ information is lost. It is not a stable homotopy category [\ref{dis:classical-stable-hotop}]; in particular, $\Sigma$ and $\Omega$ are not homotopy inverses to each other. Hence, its homotopy category has loss of spatial information, and does not possess a triangulated structure.
  \end{proof}
\end{olemma}

\begin{odiscussion}[Motivation for the machinery of \8-operads]
  The motivation for the machinery of \8-operads arises from the fact that they built out of simplicial sets, the homotopy category of whose spectra possess a triangulated structure [\ref{def:triangulated-cat}] (which endow them with spatial information), and whose spectra are homotopy categories of stable \8-categories [\ref{dis:classical-stable-hotop}, \ref{def:Sp}].
\end{odiscussion}

\begin{theorem}[Calaque-Campos-Nuiten]
  Given a formal moduli problem \FMP, the main result of \cite{Calaque19} can be stated as the equivalence of the \8-categories

  \begin{align*}
    \FMPP & \simeq \Alg_{\Pz}   \\
    F     & \mapsto \Tgt{F}[-1] \\
  \end{align*}

  where $\Tgt{F}$ is the tangent complex of $F$.
\end{theorem}

\subsection{Relationship to Lie algebras}
\begin{discussion}[Formulating the deformation problem]
  In order to prove the main result, we start with the adjunction

  \begin{equation*}
    \begin{tikzcd}
      \Dfrak : \Alg_\P \arrow[r, bend left] \arrow[r, phantom, "\bot" description] & \op{\Alg}_{\Dfrak(\P)} : \Dfrak' \arrow[l, bend left]
    \end{tikzcd}
  \end{equation*}

  where $\Dfrak(\P)$ is the \K-linear dual of the bar construction

  \begin{align*}
    \Dfrak_\K(\P) := \B_\K(\P)^\vee
  \end{align*}

  Moreover, the Koszul dual property of $\P$ asserts that there are weak equivalences (using square brackets to denote degree shift):

  \begin{equation*}
    \Dfrak(\P) \tilde{\longrightarrow} \Pz[1]
  \end{equation*}

  The above adjunction is never an equivalence, because both $\Dfrak$ and $\Dfrak'$ send an algebra to its module of derivations with coefficients in \K. To prove the main result, we impose conditions on $\Dfrak$, so that it satisfies the axioms of a formal moduli problem, and obtain an equivalence between $\Dfrak(\P)$ and \FMPP.
\end{discussion}

\begin{theorem}[Deformation problem in the operadic setting]
  Let $\P$ be an \K-operad. There is then an equivalence of \8-categories

  \begin{align*}
    \MC: \Alg_{\Dfrak(\P)} & \tilde{\longrightarrow} \FMPP            \\
    \g                     & \mapsto \Map_{\Dfrak(\P)}(\Dfrak(-), \g) \\
  \end{align*}

  if the following conditions are satisfied.

  \begin{enumerate}
    \item[(A)] For every Artin \P-algebra \A, the unit map $\A \rightarrow \Dfrak\Dfrak'(\A)$ is an equivalence.
    \item[(B)] For every trivial algebra $\K[n]$ generated by a single element of degree $n \geq 0$, the \Dfrak(\P)-algebra $\Dfrak(\K[n])$ is freely generated by $\K[n]^\vee$.
    \item[(C)] The functor $\Dfrak$ sends every pullback square of Artin algebras

      \begin{equation*}
        \begin{tikzcd}
          \A' \arrow[r]\arrow[d] & 0 \arrow[d] \\
          \A \arrow[r] & \K[n]
        \end{tikzcd}
      \end{equation*}

      with $n \geq 1$ to a pushout square of \Dfrak(\P)-algebras.
  \end{enumerate}
\end{theorem}

When $\P$ satisfies suitable finite-dimensionality conditions, the formal moduli problem $\MC_\g$ can be described in terms of Maurer-Cartan simplicial sets of Lie algebras.

\begin{discussion}[The degenerate case of $\MC_\g$\label{dis:degenerate-mcg}]
  The \K-operadic version of the result says that, for a suitable augmented \K-operad \P, every algebra $\g$ over the dual operad \Dfrak(\P) determines the formal moduli problem

  \begin{equation*}
    \MC_\g: \ArtP \rightarrow \S
  \end{equation*}

  where

  \begin{equation*}
    \MC_\g(\A) = \Map(\Dfrak(\A), \g)
  \end{equation*}

  This can be understood as Maurer-Cartan elements of nilpotent \Linf-algebras.

  To recover the classical result, set $\P = \Com$ and $\g = \Lie$ to obtain

  \begin{equation*}
    \MC_\g(\A) = \MC(\A \otimes \g \otimes \Omega_\bullet)
  \end{equation*}
\end{discussion}

\begin{discussion}[Deformation functor and the Koszul duality]
  The minimum required condition on operads is detailed in the \emph{Splendid operads} section of \cite{Calaque19}. However, it is convenient to use the Koszul duality to construct the $\B$ algebra from the $\A$ algebra. An incarnation of the Koszul duality exists in the setting of \8-operads, which is what Lurie uses in his work.
\end{discussion}

\newpage
\appendix
\section{Filtration in Category Theory\label{sec:filtration}}
\begin{definition}[Filtered category]
  A non-empty small category $\Cscr$ is termed \emph{filtered} if

  \begin{enumerate}
    \item[(i)] For every $x, y \in \Cscr$, there exists $z \in \Cscr$, and morphisms $z \rightarrow x$ and $z \rightarrow y$.
    \item[(ii)] For every $f, g: a \rightarrow b$ in \Cscr, there exists $h: c \rightarrow b$ such that $fh = gh$.
  \end{enumerate}
\end{definition}

\begin{definition}[\ind]
  The category $\ind(\Cscr)$ of \emph{ind-objects} in $\Cscr$ is a category that has all diagrams $\F: \D \rightarrow \Cscr$ where $\D$ is a filtered category. Moreover, for two objects in \ind(\Cscr), $\F: \D \rightarrow \Cscr$ and $\G: \E \rightarrow \Cscr$,

  \begin{equation*}
    \Hom_{\ind(\Cscr)}(\F, \G) := \lim_{d \in \D}\colim_{e \in \E}\Hom_\Cscr(\F(d), \G(e))
  \end{equation*}

  where (co)limits are taken in \Set.
\end{definition}

\begin{definition}[\pro]
  The category $\pro(\Cscr)$ of \emph{pro-objects} in $\Cscr$ is a category that has all diagrams $\F: \D \rightarrow \Cscr$ where $\D$ is a cofiltered category. Moreover, for two objects in \pro(\Cscr), $\F: \D \rightarrow \Cscr$ and $\G: \E \rightarrow \Cscr$,

  \begin{equation*}
    \Hom_{\pro(\Cscr)}(\F, \G) := \lim_{e \in \E}\colim_{d \in \D}\Hom_\Cscr(\F(d), \G(e))
  \end{equation*}

  where (co)limits are taken in \Set.
\end{definition}

\begin{lemma}[Equivalence and \ind\label{lemma:equiv-ind}]
  Let $\D$ be a full subcategory of category \Cscr, such that

  \begin{enumerate}
    \item[(i)] $\Cscr$ is cocomplete.
    \item[(ii)] There exist functors

      \begin{align*}
        \delta & : \Cscr \rightarrow \ind(\D) \\
        \colim & : \ind(\D) \rightarrow \Cscr \\
      \end{align*}

      such that $\colim \delta = \id_\Cscr$.
    \item[(iii)] Every object in $\D$ is compact in \Cscr.
  \end{enumerate}

  Then, functors $\colim$ and $\delta$ induce an equivalence of categories.
\end{lemma}

\begin{example}[Relationship between $\Set$ and \fSet]
  Let $\Cscr \in \Set$, and let $\D \in \fSet$ be the full subcategory of finite sets. Every set is the colimit of its finite sets, so taking the diagram of all finite subsets yields

  \begin{equation*}
    \delta: \Set \rightarrow \ind(\fSet)
  \end{equation*}

  \ref{lemma:equiv-ind} tells us that $\Set$ is equivalent to \ind(\fSet).
\end{example}

\begin{notation}[$F_\bullet$]
  The notation $F_\bullet\A$ is used to denote the \emph{filtrations} of \A. Without additional qualifiers, a filtration of $\A$ is a sequence of morphisms

  \begin{equation*}
    \ldots \rightarrow \A_n \rightarrow \ldots \rightarrow \A_2 \rightarrow \A_1 \rightarrow \A_0 \rightarrow \A
  \end{equation*}

  or

  \begin{equation*}
    \A \rightarrow \A_0 \rightarrow \A_1 \rightarrow \A_2 \rightarrow \ldots \rightarrow \A_n \rightarrow \ldots
  \end{equation*}
\end{notation}

\begin{definition}[Filtered quasi-isomorphism\label{def:filtered-qiso}]
  A \emph{filtered quasi-isomorphism} is a quasi-isomorphism that is compatible with filtration. Concretely, given algebras $\A$ and \B, a filtered quasi-isomorphism between them is a quasi-isomorphism such that morphisms

  \begin{equation*}
    \A/F_\bullet\A \rightarrow \B/F_\bullet\B
  \end{equation*}

  are also quasi-isomorphisms.
\end{definition}

\section{Algebra}
\begin{definition}[Artin ring\label{def:art-rng}]
  A ring is said to be \emph{Artin} if it satisfies the descending chain condition on ideals. The ring of dual numbers $\K[\epsilon]/\epsilon^2$ is a classic example of an Artin ring.
\end{definition}

\begin{definition}[$\Sbb_n$]
  The \emph{symmetric group}, $\Sbb_n$, is defined as the group of automorphisms on the set $\{1, \ldots, n\}$. Members of this group, written as $[\sigma(1), \ldots, \sigma(n)]$, are sometimes referred to as \emph{permutations}.
\end{definition}

\begin{definition}[Shuffle\label{def:shuffle}]
  A shuffle, denoted $\Sh(p, q)$, is a subgroup of $\Sbb_{p + q}$ that satisfies

  \begin{align*}
    \sigma(1)     & < \ldots < \sigma(p)     \\
    \sigma(p + 1) & < \ldots < \sigma(p + q)
  \end{align*}

  For example, the three shuffles in $\Sh(1, 2)$ are $[1, 2, 3]$ (the identity shuffle), $[2, 1, 3]$, and $[2, 3, 1]$.
\end{definition}

\begin{definition}[Lie algebra]
  An ordinary Lie algebra is a vector space $\g$ equipped with a bilinear skew-symmetric map $[-, -]: \g \otimes \g \rightarrow \g$, which satisfies the \emph{Jacobi identity}.

  \begin{equation*}
    \forall x, y, z \in \g: [x, [y, z]] + [z, [x, y]] + [y, [z, x]] = 0
  \end{equation*}
\end{definition}

\begin{definition}[Maurer-Cartan element of Lie algebra\label{def:mc-elem}]
  A \emph{Maurer-Cartan element} of Lie algebra $\g$ is defined as a degree $1$ element $x \in \g$ satisfying the Maurer-Cartan equation.

  \begin{equation*}
    dx + \frac{1}{2}[x, x] = 0
  \end{equation*}

  The set of Maurer-Cartan elements of $\g$ is represented as $\MC(\g)$.
\end{definition}

\begin{definition}[Canonical filtration of Lie algebra]
  The \emph{canonical filtration} of Lie algebra $\g$ is defined as

  \begin{equation*}
    F^\mathsf{Lie}_n \g := \{x \in \g \mid \text{$x$ can be obtained as a bracketing of at least $n$ elements of $\g$}\}
  \end{equation*}

  Notice that $F^\mathsf{Lie}_n \g$ is both a sub-Lie algebra and a Lie ideal of \g.
\end{definition}

\begin{terminology}[Complete with respect to a filtration\label{term:complete-canonical}]
  A Lie algebra $\g$ is termed \emph{complete with respect to its canonical filtration} if

  \begin{equation*}
    \g \cong \lim_n \g/F^\mathsf{Lie}_n \g
  \end{equation*}
\end{terminology}

\begin{definition}[Complete Lie algebra\label{def:complete-lie}]
  A (complete) filtration of Lie algebra $\g$ is the sequence of subspaces

  \begin{equation*}
    \g := F_0\g \supseteq F_1\g \supseteq F_2\g \supseteq \ldots
  \end{equation*}

  such that

  \begin{enumerate}
    \item[(a)] $\g \cong \lim_n \g/F_n\g$.
    \item[(b)] $[F_i\g, F_j\g] = F_{i + j}\g$.
    \item[(c)] $d(F_n\g) \subseteq F_n\g$.
  \end{enumerate}

  A \emph{complete Lie algebra} is defined as the pair $(\g, F_\bullet\g)$.
\end{definition}

\begin{definition}[$\Omega_\bullet$\label{def:sullivan-alg}]
  The \emph{Sullivan algebra}~\footnote{Sullivan algebra has roots in rational homotopy theory. The book \cite{Felix12} is an especially good resource.} $\Omega_\bullet$ denotes the simplicial commutative algebra of polynomial differential forms on an $n$-simplex,

  \begin{equation*}
    \Omega_n := \frac{\K[t_0, \ldots, t_n, dt_0, \ldots, dt_n]}{\sum_{i = 1}^n t_i = 1, \sum_{i = 1}^n dt_i = 0}
  \end{equation*}

  endowed with a square-zero differential.
\end{definition}

\section{Monoids and Monads}
\begin{definition}[Monoid\label{def:monoid}]
  A monoid $\M$ is an object $\M$ in monoidal category \Cscr, endowed with a composition operation $\circ$, which along with two additional pieces of data, $\gamma: \M \circ \M \rightarrow \M$ and $\eta: \I \rightarrow \M$, makes the following diagrams commute.

  
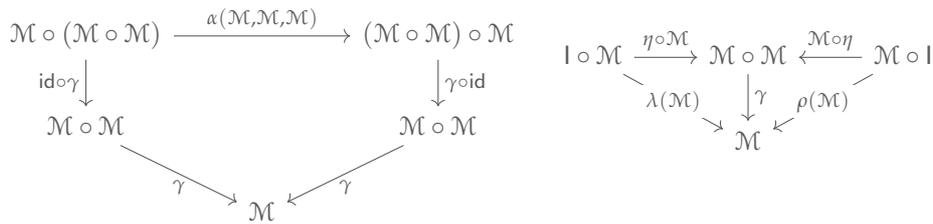
\begin{equation*}
    \begin{tikzcd}
      \M \circ (\M \circ \M) \arrow[rr, "{\alpha(\M, \M, \M)}"] \arrow[d, "{\id \circ \gamma}"{left}] & & (\M \circ \M) \circ \M \arrow[d, "{\gamma \circ \id}"{right}]\\
      \M \circ \M \arrow[dr, "{\gamma}"{below}] & & \M \circ \M \arrow[dl, "{\gamma}"{below}] \\
      & \M &
    \end{tikzcd}
    \hskip 10pt
    \begin{tikzcd}
      \I \circ \M \arrow[r, "{\eta \circ \M}"] \arrow[dr, "{\lambda(\M)}"{description}] & \M \circ \M \arrow[d, "{\gamma}"] & \M \circ \I \arrow[l, "{\M \circ \eta}"{above}] \arrow [dl, "{\rho(\M)}"{description}] \\
      & \M & \\
    \end{tikzcd}
  \end{equation*}

  where

  \begin{align*}
    \alpha(\M, \M, \M) & : \M \circ (\M \circ M) \rightarrow (\M \circ M) \circ \M \\
    \lambda(\M)        & : \I \circ \M \rightarrow \M                              \\
    \rho(\M)           & : \M \circ \I \rightarrow \M                              \\
  \end{align*}

  are natural isomorphisms in the monoidal category~\footnote{In a \emph{strict} monoidal category, this natural isomorphism is an equality}~\footnote{In a \emph{symmetric} monoidal category, $\lambda$ and $\rho$ encode the same natural isomorphism}.
\end{definition}

\begin{definition}[Monad\label{def:monad}]
  A monad $\M$ is a monoid in the strict monoidal category of endofunctors, $(\End{\Cscr}, \circ, \I_\Cscr)$.
\end{definition}

\begin{definition}[Monad algebra\label{def:monad-algebra}]
  An algebra over a monad \M, or a monad algebra, is an object $A \in \Cscr$ endowed with $\gamma_A: \M(A) \rightarrow \M$ that makes the following diagrams commute:

  \begin{equation*}
    \begin{tikzcd}
      \M \circ \M(A) \arrow[r, "{\gamma(A)}"] \arrow[d, "{\M(\gamma_A)}"{left}] & \M(A) \arrow[d, "{\gamma_A}"] \\
      \M(A) \arrow[r, "{\gamma_A}"{below}] & A \\
    \end{tikzcd}
    \hskip 10pt
    \begin{tikzcd}
      \id_\Cscr(A) \arrow [r, "{\gamma(A)}"] \arrow[dr, equal] & \M(A) \arrow[d, "\gamma_A"{right}] \\
      & A \\
    \end{tikzcd}
  \end{equation*}
\end{definition}

\newpage
\section{Trees\label{sec:trees}}
\begin{definition}[Rooted tree]
  A \emph{rooted tree} is a tree with a choice of a univalent (external) vertex.

  
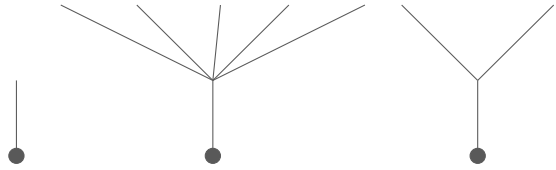
\begin{figure}[H]
    \centering
    \begin{tikzpicture}
      \draw (1, 0) -- (1, 1);
      \draw[fill] (1, 0) circle[radius=.1];
    \end{tikzpicture}
    \hskip 10pt
    \begin{tikzpicture}
      \draw (2, 0) -- (2, 1);
      \draw (2, 1) -- (0, 2);
      \draw (2, 1) -- (1, 2);
      \draw (2, 1) -- (2.1, 2);
      \draw (2, 1) -- (3, 2);
      \draw (2, 1) -- (4, 2);
      \draw[fill] (2, 0) circle[radius=.1];
    \end{tikzpicture}
    \hskip 10pt
    \begin{tikzpicture}
      \draw (1, 0) -- (1, 1);
      \draw (1, 1) -- (0, 2);
      \draw (1, 1) -- (2, 2);
      \draw[fill] (1, 0) circle[radius=.1];
    \end{tikzpicture}
  \end{figure}

  Rooted trees are denoted $\RT$.
\end{definition}

\begin{definition}[Planar tree]
  A tree is termed \emph{planar} if the ordering of the leaves is significant. Since, in our study of operads, the leaves are never labeled, this amounts to saying that the tree cannot be reflected about the $y$-axis. For example, the following two planar rooted trees are distinct.

  \begin{figure}[H]
    \centering
    \begin{tikzpicture}
      \draw (1, 0) -- (1, 1);
      \draw (1, 1) -- (0, 2);
      \draw (.5, 1.5) -- (1, 2);
      \draw (1, 1) -- (2, 2);
    \end{tikzpicture}
    \hskip 10pt
    \begin{tikzpicture}
      \draw (1, 0) -- (1, 1);
      \draw (1, 1) -- (0, 2);
      \draw (1, 1) -- (2, 2);
      \draw (1.5, 1.5) -- (1, 2);
    \end{tikzpicture}
  \end{figure}
\end{definition}

\begin{notation}[$\PT_{n, k}$]
  A \emph{planar rooted tree} is denoted $\PT_n$, where $n$ is the number of leaves. One can further classify them as $\PT_{n, k}$, where $k$ is the number of vertices. Obviously, $\PT_n = \bigcup_k \PT_{n, k}$.

  \begin{figure}[H]
    \centering
    \begin{tikzpicture}
      \node at (-1, 1) {$\PT_3 := \Biggl\{$};
      \draw (1, 0) -- (1, 1);
      \draw (1, 1) -- (0, 2);
      \draw (1, 1) -- (1.2, 2);
      \draw (1, 1) -- (2, 2);
      \node[font=\huge] at (2.5, 1) {,};
    \end{tikzpicture}
    \begin{tikzpicture}
      \draw (1, 0) -- (1, 1);
      \draw (1, 1) -- (0, 2);
      \draw (.5, 1.5) -- (1, 2);
      \draw (1, 1) -- (2, 2);
      \node[font=\huge] at (2.5, 1) {,};
    \end{tikzpicture}
    \begin{tikzpicture}
      \draw (1, 0) -- (1, 1);
      \draw (1, 1) -- (0, 2);
      \draw (1, 1) -- (2, 2);
      \draw (1.5, 1.5) -- (1, 2);
      \node[font=\huge] at (2.5, 1) {,};
    \end{tikzpicture}
    \begin{tikzpicture}
      \draw (1, 0) -- (1, 1);
      \draw (1, 1) -- (0, 2);
      \draw (.5, 1.5) -- (0.98, 2);
      \draw (1, 1) -- (2, 2);
      \draw (1.5, 1.5) -- (1.02, 2);
      \node at (2.5, 1) {$\Biggr\}$};
    \end{tikzpicture}
  \end{figure}
\end{notation}

\begin{terminology}[$k$-corolla]
  $\PT_{n, k}$ is termed the $k$-\emph{corolla} of $\PT_n$.
\end{terminology}

\begin{notation}[|v|]
  $|v|$ is used to denote the number of inputs to the vertex $v$.
\end{notation}

\begin{notation}[\vert]
  Given tree $t$, $\vert(t)$ is used to denote the set of vertices of the tree.
\end{notation}

\section{\texorpdfstring{\8}{∞}-categories}
\begin{definition}[Artin algebra\label{def:art-alg}]
  An \Einf-algebra, $\A$ will be said to be \emph{Artin} if it is connective, and $\pi_* \A$ is finite-dimensional.
\end{definition}

\begin{definition}[Triangulated category\label{def:triangulated-cat}]
  A triangulated category consists of the following pieces of data:

  \begin{enumerate}
    \item[(1)] Additive category \D.
    \item[(2)] A translation functor $\D \rightarrow \D$, which is an equivalence of categories. We denote this functor as $X \mapsto X[1]$.
    \item[(3)] A collection of distinguished triangles $X \rightarrow Y \rightarrow Z \rightarrow X[1]$.
  \end{enumerate}

  These data are required to satisfy the following axioms:

  \begin{enumerate}
    \item[(TR1)] (a) Every morphism $f: X \rightarrow Y$ in $\D$ can be extended to a distinguished triangle in \D.
      (b) The collection of distinguished triangles is stable under isomorphism.
      (c) Given an object $X \in \D$, the diagram $X \rightarrow X \rightarrow 0 \rightarrow X[1]$ is a distinguished triangle.
    \item[(TR2)] A diagram
      \begin{equation*}
        X \rightarrow Y \rightarrow Z \rightarrow X[1]
      \end{equation*}
      is a distinguished triangle if and only if the rotated diagram is a distinguished triangle:
      \begin{equation*}
        Y \rightarrow Z \rightarrow X[1] \rightarrow Y[1]
      \end{equation*}
    \item[(TR3)] Given a commutative diagram
      \begin{equation*}
        \begin{tikzcd}
          X \arrow[r]\arrow[d] & Y \arrow[r]\arrow[d] & Z\arrow[r]\arrow[d, dashed] & X[1] \arrow[d] \\
          X' \arrow[r] & Y' \arrow[r] & Z' \arrow[r] & X'[1] \\
        \end{tikzcd}
      \end{equation*}
      in which both horizontal rows are distinguished triangles, there exists a dotted arrow rendering the entire diagram commutative.

    \item[(TR4)] Suppose given three distinguished triangles:
      \begin{equation*}
        \begin{tikzcd}
          X \arrow[r] & Y \arrow[r] & Y/X \arrow[r] & X[1] \\
          Y \arrow[r] & Z \arrow[r] & Z/Y \arrow[r] & Y[1] \\
          X \arrow[r] & Z \arrow[r] & Z/X \arrow[r] & X[1] \\
        \end{tikzcd}
      \end{equation*}
      There exists a fourth distinguished triangle
      \begin{equation*}
        \begin{tikzcd}
          Y/X \arrow[r] & Z/X \arrow[r] & Z/Y \arrow[r] & Y/X[1]
        \end{tikzcd}
      \end{equation*}
      such that the diagram
      \begin{equation*}
        \begin{tikzcd}
          X \arrow[rr]\arrow[dr] && Z \arrow[dr]\arrow[rr] && Z/Y \arrow[rr]\arrow[dr] && Y/X[1] \\
          & Y \arrow[ur]\arrow[dr] && Z/X \arrow[ur]\arrow[dr] && Y[1]\arrow[ur] & \\
          && Y/X \arrow[rr]\arrow[ur] && X[1] \arrow[ur] && \\
        \end{tikzcd}
      \end{equation*}
  \end{enumerate}
  The homotopy category of a stable \8-category \Cscr, $h\Cscr$ possesses a triangulated structure.
\end{definition}

\begin{intuition}[Presentable \8-category\label{int:pres-inf}]
  A presentable \8-category is an \emph{accessible} \8-category, that admits all small colimits. An \8-category is termed accessible if it has a good supply of filtered colimits and compact objects. What this concretely means for us are the following properties.

  \begin{enumerate}
    \item[(a)] Accessible \8-categories are stable under homotopy fiber products.
    \item[(b)] Let $\F : \Cscr \rightarrow \op{\S}$ be a functor. If $\Cscr$ is presentable, it implies that $\F$ preserves colimits.
  \end{enumerate}
\end{intuition}

\begin{definition}[\8-category of spaces\label{def:inf-space}]
  Let $\Kan$ denote the full subcategory of $\SSet$ spanned by Kan complexes. Then, the simplicial nerve of $\Kan$ is termed as the \8-category of spaces, and denoted \S.
\end{definition}

\begin{notation}[\Sfin\label{not:Sfin}]
  Let $\S_*$ denote the \8-category of pointed spaces in \S. That is, $\S_*$ is the full subcategory of $\Hom(\Delta^1, \S)$ spanned by those morphisms $X \rightarrow Y$ for which $X$ is a final object of \S. Let $\S^\mathsf{fin}$ denote the smallest subcategory of $\S$ which contains the final object and is stable under finite colimits. Then $\Sfin \subseteq \S_*$ is the \8-category of pointed spaces in $\S^\mathsf{fin}$.
\end{notation}

\begin{intuition}[Sifted colimits\label{int:siftlim}]
  In ordinary category theory, sifted colimits can ``almost'' be characterized as combinations of filtered colimits and reflexive coequalizers. In \8-categories, of particular consequence to us, is the property that colimits of simplicial objects in an \8-category are sifted; this can otherwise be stated as ``simplicial \8-colimits are sifted''.
\end{intuition}

\begin{intuition}[Barr-Beck theorem\label{int:barrbeck}]
  Given functor $\F : \Cscr \rightarrow \D$, which admits a left adjoint $\G$, the \8-categorical analogue of the Barr-Beck theorem tells us that $\G \circ \F = \Hom(\Cscr, \Cscr)$ can be promoted to a monad $T \in \mathsf{Alg}(\End{\Cscr})$ on \Cscr, and that $\G$ factors through $\G' : \D \rightarrow \mathsf{LMod}_T(\Cscr)$ of left $T$-module objects of \Cscr.
\end{intuition}

\begin{terminology}[Triangle]
  Let $\Cscr$ be a pointed \8-category. A \emph{triangle in $\Cscr$} is a diagram in $\Simplex{1} \times \Simplex{1} \rightarrow \Cscr$ depicted as

  \begin{equation*}
    \begin{tikzcd}
      X \arrow[r] \arrow[d] & Y \arrow[d] \\
      0 \arrow[r] & Z \\
    \end{tikzcd}
  \end{equation*}
\end{terminology}

\begin{definition}[Stable \8-category\label{def:stable-infcat}]
  An \8-category $\Cscr$ is termed \emph{stable} if it satisfies the following conditions.

  \begin{enumerate}
    \item[(a)] There exists a zero object $0 \in \Cscr$.
    \item[(b)] Every morphism in $\Cscr$ admits a kernel and cokernel.
    \item[(c)] Every triangle in $\Cscr$ is exact if and only if it is coexact.
  \end{enumerate}
\end{definition}

\begin{theorem}[Dold-Kan correspondence\label{thm:dold-kan}]
  The Dold-Kan correspondence establishes an equivalence between simplicial sets and (homologically) non-negatively graded chain complexes of an Abelian category. The \8-categorical analogue states that, if $\Cscr$ is a stable \8-category, then the \8-categories $\Hom(N(\Z_{\geq 0}), \Cscr)$ and $\Hom(N(\op{\Delta}), \Cscr)$ are equivalent. Intuitively, it means that, to every simplicial object $X_\bullet$ of stable \8-category \Cscr, we associate the filtered object

  \begin{equation*}
    D(0) \rightarrow D(1) \rightarrow D(2) \rightarrow \ldots
  \end{equation*}

  where $D(k)$ is the $k$-skeleton of $X_\bullet$. In particular, $\colim D(j)$ can be identified with geometric realizations of $X_\bullet$.
\end{theorem}

\begin{discussion}[Background from classical stable homotopy theory\label{dis:classical-stable-hotop}]
  The study of \emph{stable homotopy theory}~\cite{Strickland20} arose from multiple pain points, but we will place emphasis on one, in particular: the parallel between homology and homotopy is incomplete in the following sense. In homology, the suspension functor induces the following isomorphism.

  \begin{align*}
    \tilde{H}_{n + k} \Sigma^k X \cong \tilde{H}_n X
  \end{align*}

  Now, the Freudenthal suspension theorem says that

  \begin{align*}
    \pi_{n + k} \Sigma^k X \cong \pi_{n + k + 1} \Sigma^{k + 1} X
  \end{align*}

  given a large enough, finite, $n$. We would achieve the parallel that we wanted if we could smash everything with $S^{-k}$, a ``negative sphere''. For this to work, we need to define $\Sigma^\8$. Let us start with the following observation. Let $A$ and $B$ be finite CW-complexes with basepoints.

  \begin{align*}
    [A, B] \xrightarrow{\Sigma} [\Sigma A, \Sigma B] \xrightarrow{\Sigma} [\Sigma^2 A, \Sigma^2 B] \xrightarrow{\Sigma} \ldots
  \end{align*}

  Except for the first two terms, this is a sequence of homomorphisms between Abelian groups. By Freudenthal's suspension theorem, after a finite number of terms, the sequence turns into a sequence of isomorphisms. We define

  \begin{align*}
    [\Sigma^\8 A, \Sigma^\8 B] := \lim_N [\Sigma^N A, \Sigma^N B]
  \end{align*}

  to capture this convenient term, and define the category of spectra~\footnote{To be precise, the Spanier-Whitehead category of finite spectra} to have objects $\Sigma^{\8 + n} A$ and morphisms

  \begin{align*}
    [\Sigma^{\8 + n} A, \Sigma^{\8 + m} B] := \lim_N [\Sigma^{N + n} A, \Sigma^{N + m} B]
  \end{align*}

  Here, $m$ and $n$ are integers. In this category, $\Sigma$ induces a self-equivalence.

  The category of spectra has a \emph{triangulated structure}.
\end{discussion}

\begin{definition}[\Sp\label{def:Sp}]
  For any \8-category that admits finite limits, consider \8-category $\Sp(\Cscr)$ of \emph{spectrum objects} of \Cscr. In the special case where $\Cscr$ is the \8-category of spaces, we recover classical stable homotopy theory. If \8-category $\Cscr$ is pointed, the spectrum objects of $\Cscr$ can be defined as the homotopy inverse limit of the tower of \8-categories:

  \begin{equation*}
    \begin{tikzcd}
      \ldots \ar[r] & \Cscr \ar[r, "{\Omega_\Cscr}"] & \Cscr \ar[r, "{\Omega_\Cscr}"] & \Cscr
    \end{tikzcd}
  \end{equation*}

  The \8-category $\Sp$ is stable and presentable, whose homotopy category recovers the classical stable homotopy category, and in particular possesses a triangulated structure (see \cite{Lurie12}, Chapter 1). $\Sp$ can also be thought of as an \8-category whose objects are cohomology theories. An alternative viewpoint is to think of $\Sp$ as the \8-category of loop spaces. $\Sp$ bears the same relationship to the \8-category of spaces as the $\Ab$ bears to $\Set$.
\end{definition}

\end{document}